\documentstyle{article}

\newcommand{\R}{{\bf R}}

\newcommand{\ppt}{{\rm pt}}

\newcommand{\Jb}{{J_\be}}
\newcommand{\th}{{\theta}}

\newcommand{\Q}{{\bf Q}}

\newcommand{\ev}{{\rm ev}}
\newcommand{\GR}{{\rm Gr\,}}

\newcommand{\Z}{{\bf Z}}
\newcommand{\C}{{\bf C}}
\newcommand{\CP}{{\bf CP}}
\newcommand{\G}{{\rm G}}
\newcommand{\p}{{\partial}}

\newcommand{\be}{{\beta}}

\newcommand{\Sig}{{\Sigma}}
\newcommand{\Om}{{\Omega}}
\newcommand{\om}{{\omega}}
\newcommand{\eps}{{\varepsilon}}

\newcommand{\Ga}{{\Gamma}}

\newcommand{\la}{{\lambda}}

\newcommand{\Jj}{{\cal J}}
\newcommand{\Ee}{{\cal E}}

\newcommand{\Hh}{{\cal H}}

\newcommand{\Mm}{{\cal M}}
\newcommand{\Pp}{{\cal P}}

\newcommand{\Tt}{{\cal T}}

\newcommand{\Hat}{\widehat}

\newcommand{\Si}{{\Sigma}}

\newcommand{\Gr}{{\rm Gr\,}}

\newcommand{\SS}{{\smallskip}}
\newcommand{\MS}{{\medskip}}

\newcommand{\NI}{{\noindent}}
\newcommand{\proof}[1]{\noindent{\bf Proof#1:\  }}

\newcommand{\QED}{\hfill$\Box$\medskip}

\newtheorem{theorem}{Theorem}[section]
\newtheorem{cor}[theorem]{Corollary}

\newtheorem{definition}[theorem]{Definition}
\newtheorem{example}[theorem]{Example}
\newtheorem{remark}[theorem]{Remark}
\newtheorem{lemma}[theorem]{Lemma}

\newtheorem{prop}[theorem]{Proposition}
\newtheorem{proposition}[theorem]{Proposition}
\newcommand{\at}{{@}}

\title{Lectures on Gromov invariants for symplectic $4$-manifolds}
\author{Dusa McDuff\thanks{Partially supported by
NSF grant DMS 9401443.} \\ State University of New York at Stony Brook \\
(dusa\at math.sunysb.edu)\\ \vspace{.1in} based on notes taken by W.Lorek}
\date{April 2, 1996}

\begin{document}

\maketitle

\begin{center} {to be published in the Proceedings of the
NATO Summer School,\\ Montreal 1995.}
\end{center}

\section*{Introduction}

Taubes's recent spectacular work setting up a correspondence between
$J$-holo\-morphic curves in symplectic $4$-manifolds and solutions of
the Seiberg-Witten equations counts $J$-holomor\-phic curves in a
somewhat new way.  The \lq\lq standard" theory concerns itself with
moduli spaces of connected curves, and gives rise to Gromov-Witten
invariants: see for example, McDuff--Salamon~\cite{JHOL},
Ruan--Tian~\cite{RT,RT2}.  However, Taubes's curves arise as zero sets
of sections and so need not be connected.  These notes are in the main
expository.  We first discuss the invariants as Taubes defined them,
and then discuss some alternatives, showing, for example, a way of
dealing with multiply-covered exceptional spheres.  We also calculate
some examples, in particular finding the Gromov invariant of the fiber
class of an elliptic surface by counting $J$-holomorphic curves,
rather than going via Seiberg--Witten theory.

For background material on symplectic manifolds and $J$-curves the
reader can consult~\cite{JHOL,INTRO} as well as the article by
F. Lalonde in this volume.  We will make passing references to
Seiberg--Witten theory, but the reader need know nothing about it to
understand most of this article.

These notes are loosely based on the lectures which I gave in
Montreal.  The treatment of Gromov invariants has been expanded, and
the material on the classification of ruled surfaces has been written
up elsewhere (in~\cite{LM,LM1}).  Here is the plan.  A more detailed
description of the contents appears at the end of Lecture 1.  I wish
to thank R. Stern and T. Parker for some helpful comments, and W.
Lorek for taking the notes, for useful discussions concerning the
material in Lecture~5, and for a careful reading of an earlier version
of this manuscript.

 \begin{description}
\item[Lecture 1] Gromov invariants: definition and examples. \dotfill 2

\item[Lecture 2]  Proof of the main structure theorem.\dotfill 7
\item[Lecture 3]  Gromov invariants: further discussion.\dotfill 13

3.1  Multiply-covered exceptional spheres\hfill 13

3.2  Components of $\Si$ \hfill 17

3.3 Examples with disconnected $K$\hfill 21

3.4 Structure of the Gromov invariants when $b_2^+ > 1$\hfill 23

\item[Lecture 4]  Spherical Gromov invariants.\dotfill 27

\item[Lecture 5]  Calculating Gromov invariants of  tori.\dotfill 32

5.1  Tori in $S^2\times T^2$ \hfill 32

5.2 Taubes's method for counting tori\hfill 33

5.3 Elliptic surfaces\hfill 35

5.4 The Gromov invariant of a fiber sum\hfill 37

\end{description}

\section{Gromov invariants: definition and examples.}
\subsection{Basic ideas}
Let $(M,\om)$ be a compact symplectic $4$-manifold with a compatible
 almost-complex structure $J$.  Given a sequence of solutions to the
 perturbed Seiberg-Witten equations for some Spin$^c$ structure $\Ga$,
 Taubes~\cite{TAU1} constructs a regular $J$-holomorphic curve $C.$
 The curve $C$ passes through $k$ generic points, where
\begin{eqnarray*}
k\; =\; k(A) & = & \frac{1}{2}(c_{1}(A) + A\cdot A). 
\end{eqnarray*}
Here $A \in H_{2}(M, \Z) $ is the homology class of $C$ and is
determined by the Spin$^c$ structure $\Ga$, and $c_1$ is the first
Chern class of the complex rank $2$ bundle $(TM,J)$.  Because it
appears as the zero section of a certain complex line bundle, the
curve $ C$ can be disconnected and can have multiply-covered
components. It might have also components which are cusp-curves or
have singularities. A natural question that arises is: what can be
said about the geometry of such a curve?

To analyse $C$  we will  parametrise it by a $J$-holomorphic
 map 
$$ 
\phi: \Sigma \to M 
$$
 from a possibly disconnected Riemann
surface $\Sigma = \coprod \Sigma_{i}$ to the manifold $M$ 
chosen so that $[\phi_*(\Si)] = A.$
 The  multiplicity of $\phi |_{\Sigma_{i}}$ on the component $\Si_i$ 
is an integer 
$m_{i}$ such that $\phi|_{\Si_i} $  may be written as a composite
$$
\phi|_{\Si_i}:\;\;\Si_i\stackrel{\psi}{\longrightarrow} \Si_i'
\stackrel{\phi'}{\longrightarrow} M,
$$ where $\psi: \Sigma_{i} \to \Sigma_{i}'$ is a branched covering map
of degree $m_{i}$ and the $J$-holomorphic map $\phi':\Si_i'\to M$ is
somewhere injective.\footnote { A ($J$-holomorphic) map $\phi:\Si\to
M$ is said to be {\em somewhere injective} if there is a point $z\in
\Si$ at which the derivative $d\phi(z)$ has maximal rank and also is
such that $\phi^{-1}(\phi(z)) = \{z\}$.  } We will assume that the
images $\phi(\Sigma_{i})$ of the different components $\Si_i$ are
distinct.  (This may be arranged by replacing several coincident
components by a single multiply-covered component.)  Thus the image
curve $C =
\phi(\Si)$ is a finite union of distinct connected curves
 $C_i = \phi_i(\Si_i)$, each with a
multiplicity $m_i\ge 1$, such that $$ A = [C] = \sum_i m_i [C_i].  $$

Two parametrizations $\phi,\Si$ and $\phi', \Si'$ are equivalent if
their images and assigned multiplicities are equal, and we denote the
equivalence class containing $\phi,\Si$ by $(\phi,\Si)$.  The pair
$(\phi, \Sigma)$ belongs to a moduli space $\Hh(A)$ which is defined
as follows.

\begin{definition}\rm
Given $A\in H_2(M,\Z)$ let $\Omega_{k}$ be a set of $k(A) =
\frac{1}{2} ( c_{1}(A) + A \cdot A )$ distinct points on $M.$ The
moduli space $\Hh(A)$ is the set of equivalence classes $(\phi,
\Sigma)$ as above such that the image $\phi(\Sigma)$ contains
$\Omega_{k}:$
\begin{eqnarray*}
\Hh(A) = \{ (\phi, \Sigma):
 \Omega_{k} \subset \phi(\Sigma) \}
\end{eqnarray*}
Moreover a pair $(\phi, \Sigma) $ will be called {\bf good} if $\phi
|_{\Sigma_{i}}$ has multiplicity $m_{i} = 1$ whenever $\phi(\Si_i)^2 =
\phi(\Si_i)\cdot \phi(\Si_i) < 0$. \end{definition}

  Note that the elements of $\Hh(A)$ are unparametrized rather then
parametrized curves.
The following theorem is due to Taubes~\cite{TAU1}.
Intuitively, it says that $k(A)$ is the maximal dimension of a
stratum in the space of all (possibly disconnected) $J$-holomorphic 
$A$-curves.

\begin{theorem}\label{goodcurves}
Suppose that $J$ is a generic $\om$-tame almost-complex structure on $M,$ and
$A\in H_{2}(M, \Z)$
a homology class.
\begin{description}
\item[(i)]
Suppose  the moduli space $\Hh(A)$ contains a good pair.
Then every pair $(\phi , \Sigma) \in \Hh (A)$ is good.  Moreover:
\begin{description}
\item[(a)] For every component $\Sigma_{i},\, $  $\phi(\Sigma_{i}) $ is an
embedded curve, disjoint from all other curves $\phi(\Sigma_{j}).$
\item[(b)] The multiplicity $m_{i}$ of $\phi|_{\Sigma_{i}}$ is one,
unless the genus
 $g(\Sig _{i}) = 1$ and $\phi(\Si_i)$ has zero self-intersection.
\item[(c)] The moduli space $\Hh(A) $ is $0$-dimensional, and
finite. 
\end{description}
\item[(ii)] If $(\phi, \Sig) \in \Hh(A),$ then the image 
$\phi(\Sig_{i}) $ of every
$\Sig_{i} $ such that 
$\phi(\Si_i)\cdot \phi(\Si_i) < 0$ is an embedded exceptional
sphere.  However its multiplicity may be $> 1$. \end{description}
\end{theorem}

The proof is deferred to the next lecture.
If the elements of $\Hh(A)$ are not good, many of the statements 
made in (i) above still hold.  The situation is fully explained in~\S3.1.

We continue here with a brief discussion of the Gromov invariant
$\Gr(A)$ and the calculation of some easy examples.  The basic idea is
that $\Gr(A)$ counts the number of elements in $\Hh(A)$ with
appropriate sign.  It is quite easy to make this precise when no
elements of $\Hh(A)$ have components which are multiply-covered.
However, multiply-covered tori are very difficult to count, and we
will postpone further discussion of this case to Lecture~5.

For the time being, let us suppose that for every $(\phi, \Si)\in
\Hh(A)$ the components $\Si_i$ are mapped with multiplicity $1$.  Then,
if $\Si_i$ has genus $g_i$, and 
its image has homology class $A_i$ and contains $k_i$
of the points of $\Om_k$, there is an evaluation map of the form
$$
\ev:\prod_i \Mm(A_i,J,g_i)\times_{G_i} (\Si_i)^{k_i} \to M^{k(A)},
$$ 
where $\Mm(A_i,J,g_i)$ is the moduli space of (connected)
$J$-holomorphic curves of genus $g_i$ in class $A_i$ and $G_i$ is an
appropriate reparametrization group.  (See Lecture 2.)  Note that
$\Mm(A_i,J,g_i)$ has a canonical orientation even in the case that it
is zero-dimensional. It follows from the proof of
Theorem~\ref{goodcurves} that the domain and range of $\ev$ have the
same dimension.  Thus there is a bijection between the subset of
$\Hh(A)$ corresponding to the given decomposition $A_i, k_i$ and the
set $\ev^{-1}(x_1,\dots,x_k)$, where $x_1,\dots,x_k$ are the points of
$\Om_k$ listed in appropriate order.  Since $\ev$ maps between
oriented manifolds, one can therefore assign a sign $\eps(\phi,\Si) =
\pm 1$ to each such element of $\Hh(A)$. Observe that this sign is
simply the product of signs which are attached to each component via
the evaluation map $$
\ev_i:  \Mm(A_i,J,g_i)\times_{G_i} (\Si_i)^{k_i} \to M^{k_i}.
$$ 
Note also that if $J$ is integrable (and regular) the evaluation
maps $\ev_i$ are holomorphic and so preserve orientation everywhere.
This means that in the K\"ahler case all the curves count with $+1$.

Here is a preliminary version of the definition of the Gromov
invariants, which is valid when there are no multiply-covered tori.
The general case is dealt with in Definitions~\ref{def:gr1}
and~\ref{def:gr2}.

\begin{definition}\label{def:gr}\rm 
Given a homology class $A \in H_2(M, \Z)$ such that $\Hh(A)$ only contains
elements with components of multiplicity $1$, 
we define the Gromov
invariant $\Gr(M,A) = \GR(A)$ by:
$$
\Gr(A)  =\sum_{\{(\phi,\Si) \in \Hh(A)\}}\eps (\phi, \Sigma)
$$
This number is independent of the choice of generic $\om$-tame $J$.
\end{definition}

 \subsection{Examples}
There are several basic examples where the Gromov invariants can be
rather easily computed.

\begin{example}\rm\label{ex:prodsph}
Let $M = S^{2} \times S^{2}$ with its standard integrable complex
structure,
and the standard product symplectic form.
 Let $A_{1} = \left[S^{2}\times
{\rm pt}\right] , $ and $A_{2} = \left[ {\rm pt} \times S^{2} \right]
.$ 
Then $k(A_{1}) = \frac{1}{2}(c_{1}(A_{1}) + A_{1} \cdot A_{1}) =
\frac{1}{2}(2 + 0) = 1,$ i.e we are counting $J$-curves in class $A_{1}$
passing through one generic point $z_0$.  There clearly is
 a unique $J$-holomorphic sphere $S$ in class $A_{1}$ passing
through the point $z_0$ and it is not hard to show that it is regular.
Moreover there cannot be another $J$-holomorphic
$A_1$-curve $S'$ through 
$z_0$ by positivity of intersections:  if there were we
would have $A_1^2 = S\cdot S' > 0$ which is absurd.  Hence $\GR(A_{1}) =
1$.

A similar argument shows that $\Gr(2A_1) = 1$.  In this case $k(2A_1)
= 2$ so that we are counting curves through $2$ generic points.
Because $J$ is a product, Theorem~\ref{goodcurves} implies that the
only elements in $\Hh(2A_1)$ are doubly covered $A_1$-spheres and
disconnected curves consisting of $2$ disjoint $A_1$-spheres.  Since
curves of the former type only go through $1$ generic point, we just
have to count the number of pairs of $A_1$-spheres through a given
pair of points.  Since these points are generic, they do not lie on
the same $A_1$-sphere and so (by positivity of intersections again)
there is exactly one such pair.

Next consider the class $A_1+A_2$.  Because $J$ is a product,
$J$-holomorphic spheres in class $A_1+A_2$ are graphs of holomorphic
maps $S^2\to S^2$ and so there is a unique such graph through $3$
generic points.  This is consistent with the fact that $k(A_{1} +
A_{2}) = \frac{1}{2}(4 + 2) = 3$, and implies that $\GR(A_{1} + A_{2})
= 1.$

In fact, it follows from Taubes' results and the wallcrossing formula
of Li-Liu that $\Gr(A) = 1$ for all nonzero $A= pA_1 + qA_2$ with
$p,q\ge 0$: see~\cite{LL1,LL2,LIU}.
\end{example} 

\begin{example}\rm\label{ex:CP2}
This time let $M = \C P^{2} $ with its standard complex structure.
 Let $L=\left[ \C P^{1} \right]$ and $A= 3L$.  Then $k(A) =
 \frac{1}{2}(9 + 9 ) = 9, $ and we are counting curves through $9$
 generic points.  The curves in class $A$ are the cubic curves --
 either (embedded) tori or rational curves with a double point or a
 cusp.  Recall that there exists a unique holomorphic torus through
 $9$ generic points, hence $\# \Hh(A) \geq 1.$ (We are in the
 integrable case here so that all signs are $+1$.)  On the other hand,
 the complex dimension of the moduli space of holomorphic curves of
 genus $g$ in class $A$ is $ ( c_{1}(A) + g -1).$ It follows that
 there is a finite number of rational curves through $8$ generic
 points, hence there are no rational curves through $9$ points in
 generic position.  Since all curves in class $A$ are either rational
 or tori we can conclude that $\Gr(A) = \# \Hh(A) = 1.$ For further
 discussion see Example~\ref{ex:grs}.
\SS

Here we have given an independent argument to show that the elements
of $\Hh(A)$ are embedded curves.  However, this is part of
Theorem~\ref{goodcurves}.  In fact, the proof of this part of the
theorem is just a more elaborate version of the argument presented
above.
\end{example} 

\begin{example}\rm\label{ex:blow}
Let $M = \C P^{2} \# \overline{{\C}P}\,\!^{2},$ and let $E $ be the
homology class of the exceptional divisor, so that $E \cdot E = -1.$ 
For any $J$ there exists a unique  $J$-holomorphic representative
$C_{E}$ of $E.$ 
Consider now $A = L + E. $ There is no connected
 $J$-holomorphic curve $C$ in that class.
 For if $\left[C\right] = L + E, $ then the intersection index
 $C \cdot C_{E} = -1 $ so that (by positivity of 
intersections) $C_{E}$ is a component
of $C.$ Now
 $k(A) = \frac{1}{2}(c_{1}(A) + A\cdot A) = \frac{1}{2}(4 + 1 - 1) =
2. $  It is easily seen that there is only one curve in class $A$
through  two generic points: it has two components, $C_{E}$ and an $L$-sphere
$\Sigma $.  The latter contains the two points of $\Om_2$ and is disjoint from
$C_E$.

Similarly, the class $L + 2E$ is represented by the disjoint union of
a sphere in class $L$ through $2$ generic points and a double cover of
$C_E$.  Thus $\Gr(L + 2E)$ should be $1$.  But this element
$(\phi,\Si)$ is not good, and so it does not appear in $\Hh(A)$.
Instead, observe that $k(L+2E) = 1$, so that there is a whole family
of curves through $k(A)$ generic points.  Thus part~(i)(c) of
Theorem~\ref{goodcurves} fails.  In fact there also is an isolated
representative of $A$ consisting of one $L-E$ curve together with a
triple cover of the $E$ curve, but now the different components
intersect so this should not contribute to $\Gr(A)$.
\MS

  An internally consistent definition of the Gromov invariants for
classes whose representation involves multiply-covered exceptional
spheres is presented in Lecture~3 below.  We will see that it suffices
to alter the definition of $k(A)$.  At this writing it is not clear
whether this definition is appropriate in the context of Taubes'
identification of the Gromov invariants with the Seiberg--Witten
invariants.  However, in the above example, we know that the
Seiberg--Witten invariant of the class $L + 2E$ is $1$, and so the
evidence points to it being the correct definition.  \end{example}

\subsection{ Further Contents}

To finish, we briefly describe the contents of the remaining lectures.

\NI$\bullet$ Lecture~2 gives the proof of  
Theorem~\ref{goodcurves}.  The argument
is basically straightforward even though it is somewhat long: it is
yet another indication that in dimension $4$ homology determines
geometry.

\NI$\bullet$  In Lecture~3 we take up two
important questions concerning the Gromov invariants.  The first (how
to deal with multiply-covered exceptional curves) arose in
Example~\ref{ex:blow} above.  We propose a definition of a modified
invariant $\Gr'(A)$ which takes care of this problem.  The second
question also appeared there, albeit indirectly.  It is the question
of how one knows that one has found all the elements of $\Hh(A)$.  The
decomposition of each $(\phi,\Si)\in \Hh(A)$ into its components
$(\phi|_{\Si_i},\Si_i)$ gives rise to a corresponding decomposition
$A= \sum_jB_j$.  (If there are no toral components, the set of $B_j$
is simply the set of homology classes represented by the components of
$\Si$.)  We discuss cases in which only one such decomposition
occurs. To what extent this is true in general is an open problem.

\NI$\bullet$  In Lecture~4 we define an analog $\Gr_s(A)$  of
$\Gr(A)$ which only counts spheres, and discuss its relation to $\Gr(A)$.  

\NI$\bullet$  In
Lecture~5 we discuss the calculation of $\Gr(A)$ in the case when  $A$ is
represented by tori.  We also give  examples to show 
why problems arise when counting
multiply-covered tori, and outline the method of counting them that Taubes
developed in~\cite{TAUTOR}.
Finally, using a $J$-holomorphic
analog of Gompf summing, we calculate
$\Gr(A)$ for the fiber class in an elliptic surface.

\section{Proof of the main structure theorem.}

This lecture is devoted to a proof of Theorem~\ref{goodcurves}.

We will work in the following set-up: 
\begin{description}
\item[(1)] $\Sig_{g} $  denotes a connected a $2$-dimensional
  manifold of genus
$g.$
\item[(2)] $\Tt _{g} $  denotes Teichm\"uller space.  Thus 
${\dim}_{\R} \Tt _{g} = 6 g - 6$ when $g>1$, and there is a smooth mapping 
\begin{eqnarray*}
j: \Tt_{g} &\to&  {\cal J}(\Sig)\\
\tau & \to & j(\tau)
\end{eqnarray*}
where ${\cal J}(\Sig)$ denotes the space of almost-complex structures
on $\Sig.$
\item[(3)] $G_{g}$ denotes the reparametrisation group. $G_{0} = {\rm
PSL}(2, \C), $  $G_{1} $ is an extension of ${\rm SL}(2, \Z) $ by the
torus   $T^{2}.$  For $g \geq 2$  the group $G_{g}$ is the mapping
class group, isomorphic to $\pi_{0}({\rm Diff}(\Sig)).$  Thus, 
$$
\dim\,G_0 = 6,\quad\dim\,G_1 = 2,\quad \dim\,G_g = 0 \mbox{ for } g>1.
$$
Note that $G_g$ is  the full group of automorphisms for a generic element 
of Teichm\"uller space, but there is a singular set (of
complex codimension $\ge 1$) of elements that have larger
automorphism groups.
\end{description}

\begin{lemma}\label{para}
Let  $ \phi: \Sig \to M $ be a  $(J, j(\tau))$-holomorphic curve i.e  
\begin{eqnarray*}
{\rm d}\phi \circ j(\tau) & = & J \circ {\rm d}\phi
\end{eqnarray*} 
Then for every 
$\gamma \in G_{g} $ the composition $\phi\circ \gamma $ is 
$(J, j(\gamma^{-1}\circ \tau))$-holomorphic.
\end{lemma}
\proof{}  This is obvious. \QED

\begin{definition}\rm\label{gmod}
For an almost-complex structure $J$ on $M,$ and a homology class $A
\in H_2(M, \Z)$ let $\Mm(A, J, g) $ denote the space of $J$-holomorphic
curves
of genus $g$ in class $A.$ More precisely,
\begin{eqnarray*}
\Mm(A, J, g) & = & \{(\phi, \tau) \in {\rm Maps}(\Sig, M) \times \Tt_{g} :
\mbox{ the curve}\,  \phi \mbox{ is}\\
& &  \quad (j(\tau),J)\,  \mbox{-holomorphic, }
 \mbox{ somewhere injective,}\,  \mbox{and} \\
& & \qquad\mbox{represents the homology class } \,  A\}
\end{eqnarray*}
\end{definition}
We have the following basic  theorem:

\begin{theorem}\label{fund}
For a generic $\om$-tame $J$ the moduli
 space $\Mm(A, J, g)$  is an oriented manifold of (real) dimension 
$ 2 ( c_1(A) + g - 1) + {\dim }\, G_{g} $.  Further, if $\ell = \ell_{g}(A) = c_{1}(A) + g -
1$, there is a well defined evaluation mapping:
\begin{eqnarray*}
\ev: \Mm(A, J, g)\times_{G_{g}} \Sig^{\ell} & \to & M^{\ell}\\
(\phi, \tau, z_{1}, \cdots , z_{\ell}) &\mapsto & (\phi(z_{1}), \cdots
, \phi(z_{\ell}))
\end{eqnarray*}
between manifolds of equal dimension $4\ell$.
\end{theorem}

The moduli space $\Mm(A, J, g)$ is used to count $J$-holomorphic
curves.  Roughly speaking, the number of $J$-holomorphic curves
through $\ell$ generic points is equal to the degree of $\ev$.  (Note
that $\ev$ maps between manifolds of dimension $4\ell$.)  A precise
statement requires compactification of $\Mm(A, J, g)$, hence
introduction of cusp-curves.  When $M$ has dimension $4$ (or $6$) the
set of points of $M$ which lie on $A$-cusp-curves always has
codimension $2$, and it follows by a standard argument that, except
possibly in the case $A\cdot A = 0, g = 1$, the map $\ev$ represents a
homology class.  (It is a pseudocycle in the language of \cite{JHOL}.)
For reasons of dimension, this homology class is a multiple
$p[M^\ell]$ of the fundamental class $[M^\ell]$, and we can figure out
what $p$ is by counting the points in the inverse image of any point
$(x_{1}, \cdots , x_{\ell})\in M^\ell$.  The case $A\cdot A = 0$, $A =
mB, m > 1$ and $g = 1$ must be treated separately since, although the
$A$-curves themselves are embedded tori (by definition the elements of
$\Mm(A,J,1)$ are somewhere injective), it is possible for these tori
to converge to multiply-covered tori in some class $kB$.  As Ruan
pointed out, this does not happen for generic $J$.  However, as Taubes
realised in~\cite{TAUTOR}, there are generic $1$-parameter
deformations of $J$ along which embedded tori in class $A = 2B$ are
absorbed by tori in class $B$.  Hence the number of tori in such a
class $A$ is not globally constant, although it is locally constant.
We will discuss this more in Lecture~5, contenting ourselves for now
with the following theorem.

\begin{theorem}\label{fund2}   In the situation of the previous theorem, given
 a generic set of points $(x_{1}, \cdots , x_{\ell}) \in M^{\ell}$ the
inverse image $\ev^{-1}(x_{1}, \cdots , x_{\ell}) $ is finite.  Moreover,
except possibly when $A^2=0, A = mB$ and $g = 1$, 
the number of points in this inverse
 image (counted with sign) is independent of
the choice of generic $J$.
\end{theorem}

Let now $C$ be a connected $J$-holomorphic $A$-curve of genus $g,$ and
multiplicity 1: \begin{eqnarray*}
\phi:(\Sig_{g}, j)  \to  M, \\
\phi_{*}\left[\Sig_{g}\right] = \left[ C\right] = A
\end{eqnarray*}
If  $J$ is generic and $\Mm (A, J, g)$ is non-empty then necessarily
$c_{1}(A) + g - 1 \geq 0.$  Moreover, the above theorems imply that there exist
finitely many such curves $C$ through $\ell_{g}(A) = c_{1}(A) + g - 1$ distinct
points.    We will need a version of the adjunction 
formula for such curves $C$.

\begin{proposition}\label{adj}
If a (connected) $J$-holomorphic curve $C$  has genus $g$ and is in the class $A$
then: \begin{eqnarray}
A\cdot A \geq c_{1}(A) + 2 (g -1)
\end{eqnarray}
with equality if and only if $C$ is embedded.
\end{proposition}
\proof{} (Sketch)
Suppose first that $C$ is immersed, with simple double points. Then
$c_{1}(A) = c_{1}(TM, J)(C) = c_{1}(TC)(C) + c_{1}(\nu_{C})(C) = 2
- 2 g + C\cdot C - 2 m,\,  $ where $m$ is the number of double points,
and $\nu $ is the normal bundle.
Hence:
\begin{eqnarray*}
A \cdot A = C \cdot C = c_{1}(A) + 2 (g - 1) + 2 m
\end{eqnarray*}
If $C$ is singular, use~\cite{LOCBE} to perturb $C$ to an immersed
curve with double points.  Every singularity contributes a non-zero
number of double points, and the proposition follows easily from the
immersed case.
\QED

For later reference, we recap the properties of $C$.
\begin{cor}  Let $C$ be a 
(connected) $J$-holomorphic curve for some generic $J$.
\begin{description}
\item{(i)} $C$ exists only if $\ell_{g}(A) = c_{1}(A) + g - 1 \geq 0.$
\item{(ii)}  There is only a 
finite number curves of of genus $g = g(C)$ and in the
class $[C]$ through $\ell_{g}(A)$ generic points.  In particular if
$\ell > \ell_g(A)$ there are no curves of this kind through $\ell$
generic points.  Thus we will say that $\ell_g(A)$ is the maximum
number of generic points which can lie on $C$.
\item{(iii)} The adjunction 
formula holds: $ A \cdot A \geq c_{1}(A) + 2(g -1),$ with
equality if and only if $C$ is embedded. 
\item{(iv)}  If $k(A) = \frac{1}{2}(A\cdot A +
c_{1}(A)) , $ then $k(A) \geq \ell_{g}(A), $ 
with equality if and only if $C$ is
embedded. 
\end{description}
\end{cor}
\MS

The above results hold for curves of multiplicity $1$.  
Suppose now that $C$ is a  curve with multiplicity
$m, $ but still connected.   Thus a
parametrization $(\phi, \Sig) $ of $C$ factors through a degree $1$
mapping $\phi': \Sig' \to M,$ and we define $g = g(C)$ to be the genus of 
the underlying simply-covered curve $\Sig'.$
 Further if $\left[C\right] = m B$ we
set
$$
\ell_{g,m}(mB) = \ell_{g}(B) = c_{1}(B) + (g - 1).
$$
As with the number $\ell_g(B)$ defined above, this number 
$\ell_{g,m}(mB)$ is the maximum number of generic points which lie on 
a curve such as $C$.

\begin{lemma}\label{mult}  Let $C$ be an $m$-fold cover of a $J$-curve in
class $B$ where $J$ is generic.
 \begin{description}
\item[(i)]If $B\cdot B \geq 0$ 
then $k(mB) \geq \ell_{g,m}(mB)$ with equality if
and only if either $m = 1 $ and $C$ is embedded, or $C$ is an $m$-fold
cover of an embedded torus of self-intersection $0$.
\item[(ii)]If $B\cdot B < 0$ then $C$ is a (possibly multiply-covered)
exceptional sphere.
\end{description}
\end{lemma}
\proof{}
Part (i) follows from a computation:
\begin{eqnarray*} k(mB)   & = &  \frac{1}{2}(m^{2} B\cdot B + m c_{1}(B)) \\
& \geq &
\frac{m}{2}(B\cdot B + c_{1}(B)) \\
  & \geq &m\; \ell_{g}(B)   \geq  \ell_{g}(B).
\end{eqnarray*}
If the first inequality holds with $m > 1$ we 
must have $B\cdot B = 0$, and if the
second holds we need $ \ell_{g}(B) = 0$.  
Hence $g = 1$.  Moreover, since $k(B) =
\ell_g(B)$  $C$ is an $m$-fold cover of an embedded curve.

As for the second part, we only  need  to observe
that two inequalities hold: 
\begin{eqnarray*}
 \ell_{g}(B)= c_{1}(B) + g - 1  \geq   0 \\
c_{1}(B) + 2(g -1)  \leq  B\cdot B   <  0.
\end{eqnarray*}
It follows that $g - 1 < 0 , $ hence $g = 0.$  Then $c_{1}(B) \geq 1,
$ and from the second inequality $c_{1}(B) = 1. $ Finally it follows
that $B\cdot B = -1, $ and so $C$ is an $m$-fold
cover of an (embedded) exceptional curve. 
\QED

Observe that the only time that $k(mB) $ is less than $ \ell_{g,m}(mB)$  is
when $B$ is represented by an exceptional sphere and $m > 1$.  This is why
multiply-covered exceptional spheres have 
to be treated separately.  Recall from
Lecture 1 that
 a pair $(\phi, \Sigma) $ is called \lq\lq good'' if $\phi |_{\Sigma_{i}}$
has multiplicity $m_{i} = 1$ whenever $\phi(\Si_i)\cdot \phi(\Si_i) < 0$.

\begin{cor}  The pair  $(\phi, \Sigma) $ is  good if and only if it contains no
component which is a  multiply-covered
exceptional sphere.
\end{cor}
\proof{}  This is  immediate from Lemma~\ref{mult}.\QED

 \begin{proposition}\label{2main}
Suppose that all of the elements $(\phi, \Sig) \in \Hh(A)$ are good. Then for
every curve $(\phi, \Sig)$: 
\begin{description}
\item[(i)]$\phi(\Sig_{i})$ is embedded, and disjoint from all other
components
$\phi(\Sig_{j}),$ except possibly if $\phi(\Sig_{i})$ is a torus with 
zero self-intersection. 
\item[(ii)]The multiplicity $m_i$ of 
$\phi|_{\Sig_{i}}$ is one, except possibly if
$\phi(\Sig_{i}) $ is a torus with zero self-intersection.
\end{description}
\end{proposition}
\proof{} Decompose $\{1,\dots, k\}$ into 
disjoint sets $I_p, p = 1,\dots,r,$ such that
the images
$$
C_p = \phi(\coprod_{i\in I_p} \Si_i)
$$ 
are connected and mutually disjoint.  If $A_i = \phi_*(\Si_i)$ and $A_p' =
\sum_{i\in I_p} A_i$ then, because $A_p'\cdot A_q'=0$, we have
$$
k(A) = k(A_{1}') + \cdots + k(A_{r}').
$$
Further, because $A_i\cdot A_j \ge 0$, for each $p$ we have
$$
k(A_p') \ge \sum_{i\in I_p} k(A_i)
$$
with equality if and only if $I_p$ has cardinality $1$.

Now let us look at the numbers $\ell_i$.  As above, $\ell_i$ should be
the maximum number of generic points on a curve of type $\phi(\Si)$.
Thus, if $\phi|_{\Si_i}$ is an $m_i$-fold cover of a curve in class
$B_i$ and of genus $g_i$, we set $$
\ell_i = \ell_{g_i,m_i}(mB_i) = \ell_{g_i}(B_i) = c_1(B_i) + g_i - 1.
$$
Then, because $\phi(\Si)$ goes through $k(A)$ generic points by assumption, we
must have $\sum_i \ell_i \ge k(A)$.  On the other hand, by Lemma~\ref{mult},
we know that $k(A_i) \ge \ell_i$ for all $i$.  Hence
\begin{eqnarray*}
k(A) = \sum_p k(A_{p}')  \geq \sum_i k(A_i) \ge \sum_i \ell_i \geq k(A).
\end{eqnarray*}
Therefore we must have equality everywhere.  In particular, each $I_p$
has cardinality $1$ which implies that all the curves $\phi(\Si_i)$
are disjoint, and $k(A_i) = \ell_i$ for all $i$.  The result now
follows immediately from Lemma~\ref{mult}. \QED
   
To complete the proof of Theorem~\ref{goodcurves} we need the following lemma.
We write $\Ee$ for the set of classes in $H_2(X)$ which are represented by 
exceptional spheres.  

\begin{lemma}\label{le:excep} The following statements
are equivalent.
\begin{description}
\item[(i)]  One element $(\phi, \Sig) \in \Hh(A) $ is good.
\item[(ii)]  Every element in $\Hh(A)$ is good.
\item[(iii)]  $E\cdot A \ge -1$ for every $E\in \Ee$.
\end{description}
\end{lemma}
\proof{}       
We will show that $(i)\Longrightarrow (iii)\Longrightarrow (ii)$.
Suppose first that there is a pair 
$(\Hat{\phi}, \Hat{\Si})  \in\Hh(A)$
which has  no components that are multiply-covered exceptional spheres.  
Then, if $\Hat{A}_j$ are the homology classes of the components of $\Hat\phi(
\Hat\Si)$, for each $E\in \Ee$ we have $E\cdot\Hat {A}_j \ge 0$ unless
$E = \Hat{A}_j$.  Hence
 there is at most one $j$ for which $E\cdot\Hat {A}_j < 0$, and for this
$j$ we have $E\cdot\Hat{A}_j = -1$. Therefore
$E\cdot A \ge -1$ for all $E\in \Ee$.  Thus $(i)\Longrightarrow (iii)$.

We prove that $(iii)\Longrightarrow (ii)$ by contradiction.
Therefore, let us suppose that $(\phi,\Si)$ does contain components
which are multiply-covered exceptional spheres.  By reordering the
components, we may suppose that these components are $\Si_i, i =
1,\dots, s,$ and that they have multiplicities $m_i>1$.  Note that
they occur in distinct classes $E_i\in \Ee$, because, by assumption,
all components of $\Si$ have distinct images under $\phi$ and, by
positivity of intersections, there is a unique $J$-holomorphic
representative of each class in $\Ee$.  Thus we may write $$ A =
\sum_{i = 1}^s m_iE_i + B, $$ where $B$ is represented by the pair
$(\phi, \Si' = \coprod_{i>s}\Si )$.  By construction, no component of
$B$ is a multiply-covered exceptional sphere or an $E_i$-curve.
Further, all the $k(A)$ generic points on $(\phi,\Si)$ must lie on the
the $B$-curve.  By the previous theorem (which applies because
$(\phi,\Si')$ is good), this implies that $k(B)\ge k(A)$.  Moreover,
we must have $E_i\cdot B \ge 0$ for all $i\le s$. For, if $E_i\cdot B
< 0$ it follows from positivity of intersections that every
representative of $B$ includes an $E_i$-curve of multiplicity at least
$1$, contradicting the definition of $B$.

Our hypothesis on $E\cdot A$ implies
in particular, that for  $i= 1,\dots,s$, 
$$
E_i\cdot A = -m_i + \sum_{j\ne i} E_i\cdot m_jE_j + E_i\cdot B \ge -1,
$$
so that
$$
 \sum_{j\ne i} E_i\cdot m_jE_j + E_i\cdot B \ge -1 + m_i.
$$
Therefore, since $E_i\cdot B \ge 0$,
\begin{eqnarray*}
2k(A) & = & c_1(\sum_i m_i E_i + B) + \sum_i (m_iE_i)^2 + B^2  \\
& & \qquad\qquad + 
\sum_im_i\left(\sum_{j\ne i} E_i\cdot m_jE_j + E_i\cdot B\right) + \sum_i
m_iE_i\cdot B\\ & \ge & \sum_i (m_i - m_i^2) + 2k(B) + \sum_i m_i(-1 + m_i) \\
& = & 2k(B).
\end{eqnarray*}
Therefore, we must have equality everywhere.  So, for each $i$, $$
A\cdot E_i = -1,\quad B\cdot E_i = 0, $$ and $$ -m_i + \sum_{j\ne i}
E_i\cdot m_jE_j = -1.  $$ If $E_i\cdot E_j\ge 1$ for some $i\ne j$, we
may suppose (by interchanging $i,j$ if necessary), that $m_j\ge m_i$.
But then $$ -m_i + \sum_{j\ne i} E_i\cdot m_jE_j \ge -m_i + m_j \ge 0
$$ 
which is impossible.  Therefore $E_i\cdot E_j = 0$ for all $i\ne
j$, which implies that $m_i = -1$ for all $i$, again contradicting our
choice of $m_i$.  Thus the lemma must hold.\QED

\section{Gromov invariants: further discussion}

We first show how to take into account multiply-covered exceptional
curves: see Definition~\ref{def:gr1}.  Next we discuss conditions
under which the surface $\Si$ in $(\phi,\Si)\in \Hh(A)$ is connected,
and give some examples (in minimal manifolds) where it is not.
Finally, we discuss the question of the uniqueness of the
decomposition of $\Si$ into its components.

 \subsection{Multiply-covered exceptional spheres}\label{ss:gr'}

We saw in Example~\ref{ex:blow} that Taubes's definition does not give the
expected answer when the class $A$ is represented by a curve which has a
multiply-covered exceptional sphere as one component.   Moreover
Theorem~\ref{goodcurves} fails in this case.  By 
Lemma~\ref{le:excep}, 
this happens if and only if 
$A\cdot E < -1$ for the class $E$ of some exceptional
sphere.  In fact, Taubes shows
in~\cite{TAU1} that this problem arises only for manifolds with $b_2^+ =
1$
\footnote
{
Recall that $b_2^+$ is the maximum dimension of a subspace of $H^2(M,\Q)$ on
which the quadratic form $a\cdot b= \langle a\cup b,[M]\rangle$ is positive
definite.
}
since otherwise $\Gr(A) = 0$ when some $E\cdot A < -1$. 
Nevertheless, it is worth attempting a better definition.

 It is not hard to deal with the problem.  The solution is to redefine
the number $k(A)$. (This amounts to looking at a different stratum of
the moduli space of all $J$-curves in class $A$.)  Before, we set $$
k(A) = \frac 12 (c_1(A) + A\cdot A).  $$ Now we set $$ k'(A) = \frac
12 \left(c_1(A) + A\cdot A + \sum_{E\in \Ee} (m_E(A)^2 -
m_E(A))\right), $$ 
where $\Ee $ is the set of classes $E$ which are
represented by exceptional spheres and where $$ m_E(A) = \max(-A\cdot
E, 0).  $$ (Think of $m_E(A)$ as the algebraic multiplicity of $E$ in
$A$.) We will look at the set $\Hh'(A)$ of pairs $(\phi,\Si)$ which
are defined as before, except now we require that $\phi(\Si)$ meets a
set $\Om'$ of $k'(A)$ generic points.

\begin{prop}\label{goodc2}
Suppose that $J$ is a generic almost-complex structure on $M,$ and
$A\in H_{2}(M, \Z)$
a homology class.
Then for any pair $(\phi , \Sigma) \in \Hh'_{J} (A):$
\begin{description}
\item[(a)] For every component $\Sigma_{i},\, $  $\phi(\Sigma_{i}) $ is an
embedded curve, disjoint from all other curves $\phi(\Sigma_{j}).$
\item[(b)] The multiplicity $m_{i}$ of $\phi|_{\Sigma_{i}}$ is one,
unless 
 $\phi(\Sig _{i})$ is a torus of zero self-intersection or an
exceptional sphere.  
\item[(c)] The moduli space $\Hh'(A) $ is $0$-dimensional,
and finite. 
\end{description}
\end{prop}

\proof{} 
As in Lemma~\ref{le:excep}, write
$$
A = \sum_{i =
1}^s m_iE_i + \sum_{j = 1}^\ell k_jF_j + B,
$$
 where the $E_i, F_j$ are the classes of the exceptional spheres in the image
$\phi(\Si)$ with multiplicities $m_i, k_j\ge 2$ and where $B$ is good. 
The $E_i$ are chosen so that
$$
 A\cdot E_i = - m_{E_i}(A) = -n_i < 0,\quad i = 1,\dots, s,
$$
and the $F_j$ are chosen so that
$$
A\cdot F_j \ge 0.
$$
Moreover, we choose the $m_i, k_j$ as large as possible so that $B$ contains no
components in the classes $E_i, F_j$.  Hence
$$
B\cdot E_i \ge 0,\qquad B\cdot F_j \ge 0.
$$

 We aim to show that $\ell = 0$ (i.e. there are no classes $F_j$),
that $m_i = n_i$ for all $i$ and that $k'(A) = k(B)$.  This will
easily imply that the $E_i$ are mutually disjoint and also disjoint
from $B$.  Since $B$ is good, the result will now follow from
Theorem~\ref{goodcurves}.

Observe  that, by definition, $\phi(\Si)$ goes
through $k'(A)$ generic points.  These must lie on the $B$ curve since
exceptional spheres do not move.  Hence $k(B)\ge k'(A)$.
To prove the converse, note first that for each $i\le s$
\begin{eqnarray}\label{eq:1}
E_i\cdot A = -m_i + \sum_{i'\ne i} E_i\cdot m_{i'}E_{i'} +\sum_j 
E_i\cdot k_jF_j + E_i\cdot B = -n_i
< 0. 
\end{eqnarray}
This implies that $m_i > m_{i'}$ for any $i,i' \le s$ such that $E_i\cdot E_i' \ne
0$. Hence by symmetry we must have
$$
E_i\cdot E_{i'} = 0,\quad   i,i'\le
s. 
$$
Further,
$$
\sum_{\{j:E_i\cdot F_j\ne 0\}} k_j \le m_i - n_i,
$$
so that
$$
\sum_{\{j:E_i\cdot F_j\ne 0\}} k_j^2 \le (m_i - n_i)^2.
$$
Therefore, if $L = \{j: F_j\cdot E_i\ne 0\; \mbox{ for some }\; i\}$,
\begin{eqnarray}\label{eq:2}
\sum_{j\in L} k_j^2 \le \sum_i(m_i - n_i)^2.
\end{eqnarray}
Equation~(\ref{eq:1}) also implies
$$
E_i\cdot (\sum_j k_jF_j + B) = m_i - n_i.
$$
Similarly, the fact that $A\cdot F_j \ge 0$ implies
$$
F_j\cdot (A - k_jF_j) = F_j\cdot(\sum_i m_iE_i 
+ \sum_{j'\ne j} k_{j'}F_{j'} + B) \ge
k_j,
$$
so that, when $j \not\in L$ we have
\begin{eqnarray}\label{eq:3}
F_j\cdot (A - k_jF_j) = F_j\cdot( \sum_{j'\ne j} k_{j'}F_{j'} + B) \ge
k_j. 
\end{eqnarray}
Thus, using equation~(\ref{eq:3}), we find
\begin{eqnarray*}
A^2 & = & (\sum_i m_iE_i + \sum_j k_jF_j + B)\cdot 
 (\sum_i m_iE_i + \sum_j k_jF_j + B)\\
& \ge & \sum_i (m_iE_i)^2 + \sum_j (k_jF_j)^2 + B^2
+ 2 \sum_i m_iE_i\cdot(\sum_j k_jF_j + B) \\
& & \qquad\qquad + 
 \sum_{j\not\in L} k_jF_j\cdot (\sum_{j'\ne j} k_{j'}F_{j'} + B)\\
&\ge & B^2 -\sum_i m_i^2 - \sum_j k_j^2 + 2\sum_i m_i(m_i - n_i) +
\sum_{j\not\in L} k_j^2\\
& = & B^2 + \sum_i (m_i^2 - 2 m_in_i) - \sum_{j\in L} k_j^2.
\end{eqnarray*}
Hence
\begin{eqnarray*}
2k'(A) & = & c_1(A) + A^2 + \sum_i (n_i^2 - n_i)\\
& \ge & 2k(B) + \sum_i(m_i - n_i) + \sum_i( m_i^2- 2m_in_i + n_i^2)
-\sum_{j\in L} k_j^2\\
& \ge & 2k(B) + \sum_i(m_i - n_i)\\
& \ge & 2k(B),
\end{eqnarray*}
where the penultimate inequality uses equation~(\ref{eq:2}).
But, as we observed earlier, $k(B)\ge k'(A)$.  Therefore we must have
equality everywhere.  This gives $m_i = n_i$ for all $i$, which, by
equation~(\ref{eq:1}), implies that 
$E_i\cdot F_j = E_i \cdot B = 0$ for all $i,j$.
Using equation~(\ref{eq:2}) we see also that $L = \emptyset$. Therefore,
if $B' = \sum_jk_jF_j + B$, we have that
$E_i\cdot B' = 0$ for all $i$, which easily implies that
$$
E\cdot B' \ge -1
$$
for all $E\in \Ee$.  The result now follows from Lemma~\ref{le:excep}.
\QED

We can now define modified Gromov invariants.  

\begin{definition}\label{def:gr'}\rm 
Given a homology class $A \in H_2(M, \Z)$ such that $\Hh'(A)$  contains
no multiply-covered tori,
we define the Gromov
invariant $\Gr'(A)$ by:
$$
\Gr'(A)  =\sum_{\{(\phi,\Si) \in \Hh'(A)\}}\eps (\phi, \Sigma) 
$$
Here we assign the sign $+1$ to each 
multiply-covered exceptional sphere and then define the sign $\eps (\phi,
\Sigma)$ as before.   This number $\Gr'(A) $ is independent of the choice of
generic $\om$-tame $J$. \end{definition}

\begin{lemma}\label{le:exc}  
$\Gr'(A) = \Gr(A)$ unless there is an $E\in \Ee$ such
that $E\cdot A < -1$, in which case $\Gr'(A) = \Gr(B) $ where
$$
B = A - \sum_{E:E\cdot A < - 1} (E\cdot A) E.
$$
\end{lemma}
\proof{}  By  Proposition~\ref{goodc2}, each element $(\phi,\Si)$ in
$\Hh'(A)$ consists of a good representative of the class $B$ together with a
collection of disjoint multiply-covered 
exceptional spheres, lying in the classes $E$
such that  $E\cdot A < - 1$.  The result follows immediately.\QED

\begin{remark}\rm As we shall explain in more detail in Lecture~4, somewhere
injective spheres are always assigned the 
sign $+1$.  Hence it is consistent also to
assign $+1$ to all exceptional spheres, including the multiply-covered ones.
\end{remark}

\begin{example}[Example~\ref{ex:blow} revisited]\rm  
Consider $M = \C P^2 \# \overline{{\C}P}\,\!^{2}$.  Then it
follows immediately from the above lemma that $\Gr'(L+ 2E ) = 1$.
\end{example}

\subsection{The components of $\Si$.}

 The decomposition of $(\phi,\Si)\in
\Hh(A)$ into its $\ell$ components 
$(\phi|_{\Si_i},\Si_i), i = 1,\dots,\ell,$ gives rise to
a corresponding decomposition $A= \sum_{i = 1}^\ell A_i$ where $A_i$
is the class represented by $\phi|_{\Si_i}$.  We now look at what we
can say about the $A_i$.  Are there any conditions under which $\ell =
1$?  Are the $A_i$ uniquely determined by $A$?
\MS

\subsubsection{Components of negative self-intersection}

The question of whether there are components with $\phi(\Si_i)^2 < 0$
and of how they appear is completely answered by the structure
theorems.  If we are dealing with the original invariant $\Gr(A)$ and
if $(\phi,\Si)$ is good then it follows from Theorem~\ref{goodcurves}
that the only negative components are exceptional spheres.  Moreover
$(\phi,\Si)\in \Hh(A)$ has a component which is an exceptional curve
in class $E$ if and only if $E\cdot A = -1$, and all $E$ which appear
in this way are disjoint.  Therefore, if $$ B = A - \sum_{\{E\in \Ee:
E\cdot A = -1\}} E, $$ 
there is a bijective correspondence between the
elements of $\Hh(A)$ and of $\Hh(B)$.  Further, no components of
negative self-intersection appear in $\Hh(B)$.  Therefore, we can
replace the study of the structure of elements of $\Hh(A)$ by that of
elements of $\Hh(B)$.

 Similarly, if we are dealing with 
$\Gr'(A)$,  $(\phi,\Si)$ has a component which is an $m$-fold cover of an
exceptional curve in class $E$ if and only if $E\cdot A = -m$.  Again, all $E$
which appear in this way are disjoint and all the other components of
$\phi(\Si)$ have nonnegative self-intersection.  Hence, as before, 
the structure of the components of negative self-intersection is determined by
homological information. 

Since the only difference between $\Gr(A)$ and $\Gr'(A)$ is in the negative
components, from now on we consider only $\Gr(A)$.
\MS

\subsubsection{ Components of zero self-intersection}  

These are either tori or spheres, since when the genus $g$ is $> 0$
the moduli space of embedded $J$-curves of genus $g$ and zero
self-intersection has negative dimension.  Moreover, these can give
rise to disconnected $(\phi,\Si)$.  We saw this in
Example~\ref{ex:prodsph} with spheres.  In \S\ref{ss:tor} we give a
similar example (on $T^2\times S^2$) with tori.  The next lemma shows
that it is impossible for both spheres and tori to occur.

\begin{lemma}\label{le:sph0}  Suppose that $(\phi,\Si)\in \Hh(A)$ contains a
component which is a sphere $C$ of zero self-intersection.  Then $M$ 
is a blow-up of a ruled surface with $C$ as one of the fibers.  Moreover any
other  components in $(\phi,\Si)$ of nonnegative self-intersection are
also fibers.
\end{lemma}
\proof{}  The first statement follows from the basic structure theorem
in~\cite{RR}.  It is easy to see that $[C]$ has nonempty intersection
with every other class $B$ with $B^2 \ge 0$ that could have a
$J$-holomorphic representative.  (Use the Light Cone lemma
(Lemma~\ref{le:light}) stated below.)  Hence if there are any more
components in $(\phi,\Si)$ with nonnegative self-intersection they
must also be fibers.\QED
\MS

\subsubsection{Components of positive self-intersection}  

 We will consider the cases $b_2^+ = 1$ and $b_2^+ > 1$ separately,
since they are rather different.

\subsubsection*{The case $b_2^+ = 1$}

The most relevant fact when considering the components of $K$ is the light
cone lemma.  It is useful to consider the positive cone
$$
\Pp = \{B\in H_2(M,\R): B^2 > 0\}.
$$
Since $b_2^+ = 1$ this has two components which are separated by the hyperplane
where $\om = 0$.  
The component on which $\om$ 
is positive is called the forward positive cone and is denoted by
$\Pp^+$.  Its closure is
$$
\overline{\Pp^+} = \{B\in H_2(M,\R): B^2 \ge 0, \om(B) \ge 0\}.
$$

\begin{lemma}[Light Cone lemma]\label{le:light}
Suppose that $(M,\om)$ is a symplectic $4$-manifold 
with $b_2^+ = 1$ and let $B_1,
B_2\in\overline{\Pp^+}$.  
Then $B_1\cdot B_2 > 0$ unless $B_1 = \la B_2$ and $B_1^2 = B_2^2 = 0$.
 \end{lemma}
\proof{} There is a basis $L,E_1,\dots, E_\ell$ for $H_2(M,\R)$ which is
orthogonal with respect to the intersection pairing and is such that
$L^2 = 1, E_j^2 = -1$ for all $j$.  Moreover, by changing the sign of
$L$ if necessary we may suppose that $\om(L) > 0$, i.e that $L\in
\Pp^+$.  Then the elements of $\Pp$ have the form $mL +
\sum_i\la_iE_i$ where $\sum \la_i^2 < m^2$.  Since $L\in \Pp^+$, this
element is in $\Pp^+$ exactly when $m > 0$.  Hence we may write the
$B_i$ as: $$ B_1 = mL + \sum_j\la_j E_j,\quad B_2 = nL + \sum_j
\mu_jE_j, $$ where $$ m,n > 0,\quad m^2 \ge \sum_j \la_j^2,\quad n^2
\ge \sum_j \mu_j^2.  $$ Therefore
\begin{eqnarray*}
B_1\cdot B_2  & = &  mn - \sum_j \la_j\mu_j\\
& \ge & mn - (\sum_j \la_j^2)^{\frac 12}(\sum_j\mu_j^2)^{\frac 12}\\
& \ge & 0
\end{eqnarray*}
as claimed.  Moreover, equality occurs only if all the $\la_j$ are
equal, all the $\mu_j$ are equal and if $B_1^2 = B_2^2 = 0$.  The
conclusion readily follows.
\QED

The next proposition shows that if $(\phi,\Si)$ has a component of positive
self-intersection, then this is the only one other than
exceptional curves.

 \begin{prop}\label{prop:b+1} Suppose that $b_2^+ = 1$ and consider
$(\phi,\Si)\in \Hh(A)$.  Let $\Si_i, i = 1,\dots, p$ be the components
for which $\phi(\Si_i)^2 \ge 0$.  \begin{description}
\item{(i)}  If some 
$\phi(\Si_i)^2 > 0$ then $p$ is at most $1$ and this component contains all the
$k(A)$ generic points.   
\item{(ii)}  If some 
$\phi(\Si_i)^2 = 0$ then $p$ can be $> 1$ but the classes
$[\phi(\Si_i)], i = 1,\dots, p$ differ by at most a constant factor.
Moreover, either they are all represented by spheres (in which case
$k(A) = p$ and all the classes $[\phi(\Si_i)], i = 1,\dots,p,$ are
equal) or they are all represented by tori (in which case $k(A) = 0$
and the classes $[\phi(\Si_i)], i = 1,\dots,p,$ all lie on the same
ray in $H_2$).
\item{(iii)}  If $\Si_i$ is a 
sphere  for some $i\le p$, then $M$ is a blow-up of
a rational or ruled surface. Moreover $\Gr(A) = 1$.
\end{description}
\end{prop}  \proof{}  Parts (i) and (ii) follows immediately from
Lemmas~\ref{le:sph0} and~\ref{le:light}, and the fact that components
with negative self-intersection are rigid so that they do not go
through any generic points.  Recall also that the genus of the
representing curves is determined homologically through the adjunction
formula.  The first statement in (iii) follows immediately from the
main theorem of~\cite{RR} and holds without the assumption that $b_2^+
= 1$. The second may either be proved using Seiberg--Witten theory or
by direct calculation.  See Propositions~\ref{prop:sgr1}
and~\ref{prop:sgr3} below.  \QED

\begin{remark}\rm  In fact, we have not yet defined $\Gr(A)$ in the case when 
some components in $(\phi,\Si)\in \Hh(A)$ are multiply covered tori.
The above proposition shows that in this case the components of
$(\phi,\Si)$ are either exceptional spheres or are tori whose homology
classes lie in some ray in $H_2$.  The exceptional spheres do not
affect the value of $\Gr(A)$, and so we can suppose that there are
none.  Then, we define $\Gr(A)$ to be $\Gr_0(A)$ as given in
Definition~\ref{def:gr2}.  \end{remark}

 The Gromov invariants for symplectic
manifolds with $b_2^+ = 1$ can be completely calculated thanks to the
wall-crossing formula in Seiberg--Witten theory: see Li-Liu~\cite{LL1,LL2}.  
When $H_1(M, \R) = 0$, $\Gr(A)$ is either $0$ or $1$, but if $b_1(M)\ne0$
the invariant can take different values.
 This leads to many interesting results.  For
example,  Liu showed in~\cite{LIU} that a minimal
symplectic $4$-manifold with $K^2 < 0$ is ruled.  
However, there are still several open questions about their structure: see the
survey article~\cite{MS}.

 \subsubsection*{The case $b_2^+ > 1$}

When $b_2^+ > 1$ the situation is more complicated.  For simplicity,
we will restrict attention to the minimal case.\footnote { In fact, if
$M$ is a symplectic $4$-manifold with $b_2^+ > 1$, $M$ has a unique
minimal reduction $M'$, ie there is a unique maximal set of
exceptional curves in $M$ (see~\cite{IMM}).  Moreover, there is a sum
formula which allows one to recover the Seiberg--Witten (or Gromov)
invariants of $M$ from those of $M'$.  Hence we do not lose any
information by restricting to the minimal case.}  Using
Seiberg--Witten theory, Taubes~\cite{TAU1} has proved the following
important structure theorem for Gromov invariants.  Recall that the
canonical class $K\in H_2(M)$ is the Poincar\'e dual of minus the
first Chern class of $M$, ie $$ K = - PD(c_1(TM, J)).  $$ In
particular, $2k(K) = c_1(K) + K^2 = -K^2 + K^2 = 0$.

\begin{theorem}[Taubes]\label{Kmin}  Let $M$ be a minimal symplectic manifold
with $b_2^+ > 1$.  Then \begin{description}
\item{(i)}
$\Gr(A) = 0$ except possibly if $k(A) = 0$.  
\item{(ii)}   $|\Gr(K)| = 1$. 
 \item{(iii)}  For all $A\in
H_2(M)$, $\Gr(A) = \pm \Gr(K - A)$.  
\item{(iv)}  If $K^2 = 0$ and $\Gr(A) \ne 0$
then $A^2 = 0$.  
\item{(v)} {\bf[Witten~\cite{WIT}]}  If $M$ is K\"ahler and $K^2 >
0$ then $\Gr(A) \ne 0$ only in the case  $A = 0, K$. \end{description}
\end{theorem}

I know  no  way of proving the above results just in the context of
holomorphic curves: at present one has to go via Seiberg--Witten
theory.  Note also that (ii) implies  that $K^2 \ge 0$
and $\om(K) > 0$, i.e. $K$ is in
the closure of the forward positive cone $\overline{\Pp^+}$.
This follows from
Theorem~\ref{goodcurves} on the structure of elements of $\Hh(A)$,
which states that the only components of $\phi(\Si)$
with negative self-intersection are
exceptional spheres.

\subsection{Examples with disconnected $K$}

Before going further, we look at some examples in which $K$ is
realised by a disconnected curve. The easiest example is that of
elliptic surfaces.  In this case, $K^2 = 0$ and $K$ is realised by a
disjoint union of parallel tori: see Lecture~5.  Here is another
example in which $K^2 > 0$.

\begin{example}\rm  We construct a symplectic manifold with a disconnected
representative of $K$ by the process of the Gompf sum.  Recall from
\cite{GOM}  that if $(M_i,
X_i)$ are two manifold/submanifold pairs such that the $X_i$ are
symplectically embedded surfaces of the same genus but opposite
self-intersection number, one can form their connected sum $$ M =
M_1\#_{X_1 = X_2} M_2, $$ by cutting out suitable neighborhoods of the
$X_i$ and gluing their complements together.  This is particularly
easy when the $X_i$ are tori of zero self-intersection: see~\S5.4
below.  In this case we also have $$ K_{M_i}\cdot X_i = 0,\quad i =
1,2, $$ so that $K_{M_i}$ may be represented by a cycle which is
disjoint from $X_i$.  It is then not hard to check that the canonical
class $K_M$ of $M$ is given by the formula $$ K_M = K_{M_1} + K_{M_2}
+ X_1 + X_2.  $$ (Note that this formula makes no sense when $X_i^2
\ne 0$ since none of the classes on the RHS can be identified in the
homology of the glued manifold $M$.)

As an example, consider $ T^4$ with the symplectic form $\om =
dx_1\wedge dx_2 + dx_3\wedge dx_4 + dx_1\wedge dx_3$.  Then $T^4$
contains {\it disjoint} nonparallel symplectically embedded tori
$X,Y$.  (For example, take $X = \{(x_1,x_2, 0,0)\}$ and $Y = \{(x_1,0,
x_3, 1/3)\}$.)  As in Lecture~5, let $V = V(1)$ denote the rational
elliptic surface $\C P^2$ with $9$ points blown up and fiber $F$, and
consider the triple sum $$ M = V\#_{F = X} T^4 \# _{Y = F'} V', $$
where $V'$ is another copy of $V(1)$.  Then, because $K_{T^4} = 0$ and
$K_V = -F$, the above formula shows that $$ K_M = F + F' = X + Y.  $$

To get an example with $K^2 > 0$, consider the manifold $$ M =
V(4)\#_{S(-4)= Q} \C P^2.  $$ Here $Q$ is the quadric in $\C P^2$ and
$S(-4)$ is a sphere of self-intersection $-4$ in the elliptic surface
$V(4)$.  (The manifold $V(4)$ is described in more detail in
Lecture~5.  The sphere $S(-4)$ is a section of the map $V(4)\to
\C P^1$, and $M$ is  called a rational blowdown of $V(4)$: see, for
example,~\cite{GOM}.)  The canonical class for $V(4)$ is $2F$ where
$F$ is the fiber class (represented by a torus with zero
self-intersection) and the canonical class for $\C P^2$ is, of course,
$-3L$, where $L = [\C P^1]$.  Consider the curve $C$ of genus $g_C =
2$ which is obtained by gluing a sphere in class $L$ to the fibers $F$
through the two points where $L$ meets $Q$.  Thus $C$ is made from two
copies of $T^2 - (disc)$, each with trivial normal bundle, plus a copy
of $S^2 - (2\; discs)$ which has self-intersection $+1$.  Thus $C^2 =
1$.  It is not hard to verify that $K_M = C$.  For example the
adjunction formula for $C$ works out: $$ - C^2 = -K_M\cdot C = 2 -
2g_C + C^2 = -1.  $$

To get a manifold with disconnected $K$, observe that $V(4)$ contains
many Lagrangian tori $Y$ which are disjoint from $F$.  To see this,
think of $V(4)$ as the fiber sum $V(2)\#_F V(2)$ of two copies of the
$K3$ surface $V(2)$, and realise $V(2)$ as the Kummer surface, which
is obtained from $T^4$ by identifying $(x_1,x_2,x_3,x_4)$ with
$(-x_1,-x_2,-x_3,-x_4)$ after having blown up the $16$ fixed points of
this involution.  The torus $Y = \{(x_1,0,x_3,1/3)\}$ (which is
Lagrangian for the usual symplectic form) descends to a torus in
$V(2)$ which is disjoint from a generic fiber $F$ of the projection
$V(2) \to \C P^1$ given by $(x_1,x_2,x_3,x_4)
\mapsto (x_3,x_4)$.  Hence $Y$ also embeds in $V(4) = V(2)\#_F V(2)$. 
Observe also that because $Y$ when considered as a subset of $T^4$
does not go through the $16$ fixed points of the involution, the image
of $Y$ in $V(2)$ is disjoint from the sections $S(-2)$ of
self-intersection $-2$ (which are the images of the blown-up points).
Hence we may assume that $Y$ in $V(4)$ is disjoint from the section
$S(-4)$.  Since the homology class of $Y$ in $V(4)$ is nonzero, we may
slightly perturb the symplectic form on $V(4)$ to make $Y$ a
symplectic torus.  (This is Gompf's trick: see~\cite{GOM}.)  Then we
can form the triple sum $$ M' = V(1)\#_{F = Y} V(4) \#_{S(-4)= Q} \C
P^2, $$ which has canonical class $K' = Y + K_M$ with $(K')^2 = 1$.
\end{example}

Here is another example with a disconnected $K$ which
contains no toroidal components. It was suggested to me
by Ron Stern.

\begin{example}\rm
Take two  surfaces $X,Y$ of general type which contain the Gompf
nucleus of the $K3$ surface.  (This nucleus is the union of
 a symplectic torus of square $0$  with
a symplectic sphere of square $-2$, and its regular neighborhood is just the 
trace of $0$-framed surgery on the right-handed trefoil and $-2$ surgery
on a meridonal curve.) There are plenty of such surfaces
 in, for example, complete
intersections. Now take the fiber sum 
$$
Z=X\#_{T_X = T_Y}\, Y
$$
 of these two surfaces
along the tori $T_X$ and $T_Y$ of square 0. Then 
$K_Z=K_X+K_Y+T_X+T_Y$.
The sum of the two $-2$-spheres
in the Gompf nuclei is a sphere of square $-4$
that intersects each of $T_X$ and $T_Y$ once.
  Therefore, one can form the 
connected sum $W$ of $Z$ with $\CP^2$ by identifying the 
complement of this $-4$ sphere with the complement of
the quadric surface.  (This 
is the rational blow-down process of~\cite{GOM}
 and Fintushel--Stern~\cite{FS}.) Then it is not hard to check that
$$
K_W=K_X+K_Y+C 
$$
where 
$C^2=1$: for more details see~\cite{FS}.
\end{example}

As we shall see in Lemma~\ref{le:gent} below, this phenomenon of
disconnected $K$ with $K^2 > 0$  cannot occur for minimal K\"ahler surfaces of general
type.  These manifolds  satisfy the Noether inequality $c_1^2 \ge 2\xi - 6$.
(Here $\xi$ denotes the holomorphic Euler characteristic $
\frac{1 - b_1 + b_2^+}2$.)
Study of the known examples of symplectic manifolds with $K^2 > 0$ 
has led Fintushel and Stern to
suggest that all minimal symplectic manifolds with $K$ connected must satisfy
the inequality $c_1^2 \ge \xi - 3$.

\subsection{Structure of the Gromov invariants when $b_2^+ > 1$}

In this section we show how the invariant $\Gr(A)$ is built up from a
simpler invariant which I will call $\Gr_0(A)$.  Roughly speaking,
$\Gr_0$ counts connected curves.  We will suppose that we are working
on a minimal manifold $M$ with $b_2^+ > 1$, so that the only classes
with nonzero Gromov invariants are those with $$ k(A) = \frac{-K\cdot
A + A^2}2 = 0.  $$ 
Of course, similar definitions can be made in the
case $b_2^+ = 1$.  However, the situation there is fully described in
Proposition~\ref{prop:b+1} and the remarks that follow it.

\begin{lemma}  Consider
$(\phi,\Si)\in \Hh(A)$ and let $A_i = \phi_*[\Si_i]$.
If $k(A) = 0$ and $\Gr(A) \ne 0$ then $k(A_i) = 0$ for all $i$.
Further, the genus $g_i$ of $\Si_i$ is $1 + A_i^2$.
\end{lemma}
\proof{} Observe that $A_i\cdot A_j = 0$ when $i\ne j$
since distinct components are disjoint.
Hence
$k(A) = \sum_i k(A_i)$.  Since
$k(A_i)\ge 0$ for all $i$ in order to have a nontrivial Gromov invariant,
this shows that $k(A_i) = 0$.
The last statement follows from the adjunction formula.  We
already know that all components of $\Si$ are embedded, 
and so 
$$
g_i = 1 + \frac 12(K\cdot A_i + A_i^2) = 1 + A_i^2 - k(A) = 1 + A_i^2,
$$ as claimed.\QED

 \begin{definition}\label{def:decomp}\rm  
Each element $(\phi,\Si) \in \Hh(A)$ determines a {\bf
decomposition} $D = \{B_1,\dots, B_\ell\}$ of $A$ in the following way.   If
$\Si_i$ is a component of $\Si$  of genus $\ne 1$, then the corresponding
homology class $\phi_*[\Si_i]$ is in $D$, but if $\Si_i$ has genus $1$ then we
group together all the components with homology class on the ray
$\{\la(\phi_*[\Si_i]):\la > 0\}$ into one element in $D$.  
Note that by Proposition~\ref{prop:b+1} there are no components of genus $0$.
Moreover, by the previous lemma,
any component of genus $1$ must have self-intersection $0$.

Thus the elements of $D$
are characterised by the following properties:
\smallskip

\NI
$\bullet$ $\sum_j B_j = A$ and $B_i\cdot B_j = 0, i\ne j$;

\NI
$\bullet$  If $i\ne j$ then $B_i \ne \la B_j$ for any $\la > 0$;

\NI
$\bullet$   if $B_i^2 > 0$, $B_i$ is represented by a connected and 
embedded  $J$-holomorphic submanifold;
 
\NI
$\bullet$ if $B_i^2 = 0$, $B_i$ is represented by a union of coverings
of embedded $J$-holomorphic tori whose homology classes all lie on the
ray $\{\la B_i: \la > 0\}$.
\end{definition}

For each such decomposition $D$ of $A$ we can add up (with signs) the
$J$-holomorphic representatives of $A$ with components in these
classes, getting an invariant which we will call $\Gr_D(A)$.  To be
more precise, consider the following definitions.

\begin{definition}\label{def:gr0}\rm (i)  If  $A^2 > 0$  then $\Gr_0(A)$ 
is defined to be the number of {\it connected, embedded} $J$-holomorphic
curves  of genus $1 + A^2$ in the
class $A$ counted with appropriate sign.  (This sign is determined by the
evaluation map as described in Lecture~1.)  Thus in this case the
invariant concides with the one considered by Ruan in~\cite{RUAN2}.

\NI
(ii) If $A^2 = 0$ then $\Gr_0(A)$ counts the number of representatives
of $A$ by disjoint unions of possibly multiply-covered tori with
homology classes on the ray $\{\la A:\la > 0\}$.  (This is the
invariant which is called ${\rm Qu}\,(A)$ in~\cite{TAUTOR}.)  In order
to get a number which is invariant under symplectic deformation it is
necessary to weight each component torus by a number which depends on
certain twisted Cauchy-Riemann operators in the normal bundle of the
torus.  This weighting is describes in more detail in Section~5.2.
Its possible values are $0, \pm 2,$ and $\pm 2k+1, k \ge 0$.  One of
Taubes's interesting discoveries in~\cite{TAUTOR} is that it is
impossible to get a well-defined invariant if one restricts attention
just to connected toral representatives in a fixed homology class.
\end{definition}

A class $B$ with $\Gr(B) = \Gr_0(B)$ will be called {\bf indecomposable}. 
Observe also that part (i) of
Theorem~\ref{Kmin} implies that if $\Gr_0(A) \ne 0$ then $K\cdot A
= A^2$.

We can now give a more precise definition of the Gromov invariant.

\begin{definition}\label{def:gr1}\rm   
Given a decomposition $D=  \{B_1,\dots, B_\ell\}$ of $A$  define
$$
\Gr_D(A) = \prod_{i = 1}^\ell \Gr_0(B_i),
$$
and set
$$
 \Gr(A) = \sum_D
\Gr_D(A). 
$$ 
\end{definition}

It is not known in general whether there is a unique decomposition $D$
such that $\Gr_D(K)\ne 0$.  However, this does hold for minimal
K\"ahler surfaces of general type, ie surfaces with $b_2^+ > 1$ and
$K^2 > 0$.

\begin{lemma}\label{le:gent}  If $M$ is a minimal K\"ahler
 surface of general type
 then $\Gr_D(K) = 0$ if $D\ne\{K\}$.  Hence $\Gr(K) = \Gr_0(K)$.
\end{lemma} \proof{}   One way to prove the first
statement is to use Witten's result in part (v) of Theorem~\ref{Kmin}.
 If $\Gr_D(K)\ne 0$ for some decomposition $D = \{B_1,\dots,B_\ell\}$
 then, by the very definition of $\Gr_D(K)$ we must have
 $\Gr_0(B_j)\ne0$ for all $j$.  Hence $\ell = 1$ by Witten's result.
 Another way to see this is to use the Hodge index theorem.  The
 classes $B_j$ would have to lie in $H^{1,1}(M)$ (because they can be
 represented by holomorphic curves), and also must have $B_j^2 \ge 0$
 and $K\cdot B_j \ge 0$. The result now follows from the Light Cone
 Lemma~\ref{le:light} because the intersection pairing on $H^{1,1}(M)$
 has type $1\oplus -1\oplus\dots \oplus -1$. \QED
                                                                          
As the next proposition shows we know much less about the general situation.
Here the possible presence of $J$-holomorphic tori causes extra problems.

\begin{prop}\label{prop:dec}  Let $M$ be a minimal  symplectic manifold with
$b_2^+ > 1$.   
\begin{description}
\item{(i)}  Suppose  that $\Gr_D(K) \ne 0$  for a
unique decomposition $D = \{B_1,$$\dots, B_\ell\}$.  If there is a 
 class $A$ such that $A^2 > 0$ and $\Gr_0(A) \ne 0$,  then $A$
must equal some $B_i$ and satisfy 
$$
\Gr_0(A) = \Gr(A) = \pm 1.
$$
In particular, if  $\Gr_0(T) = 0$ for all classes $T$ with $T^2 = 0$, then
 $\Gr(A) = 0$ unless $A$ is a union
of some of the $B_i$ in which case $\Gr(A) = \pm 1$.

 \item{(ii)}  If there are distinct
classes $A_1, A_2$ such that  
$$
A_1\cdot A_2> 0,\quad \Gr(A_i) \ne 0, \;\;i = 1,2,
$$
then $A_i^2 > 0$ for $i = 1,2$ and $\Gr_D(K) \ne 0$ for at least two different
decompositions.  The converse holds if $K^2 > 0$.
\end{description}
\end{prop}
\proof{}    Suppose first that $\Gr_0(A) \ne 0$.  Because $\Gr(K - A) \ne 0$
by Theorem~\ref{Kmin}
(iii),  Definition~\ref{def:gr1} implies that there is a
decomposition $D'= \{B_1',\dots, B_j'\}$ of $K - A$ such that 
$\Gr_{D'}(K-A)\ne 0$.  
Observe that because $2 k(A) = -K\cdot A + A^2 = 0$, 
$$
A\cdot(K-A) = 0.
$$
Therefore the union $\{A\}\cup D'$ is a 
decomposition of $K$, except possibly if
$A^2= 0$.  In the latter case there may be a
component  of $D'$ in the ray $\la A$, and if there is 
it must be amalgamated with
$A$.  We will suppose that done, if necessary, and call the resulting
decomposition $\overline D$.     It follows immediately that 
$$ 
\Gr_{\overline D}K \ne 0,
$$
unless we had to amalgamate some $B_i'= \la A$ with $A$ and it happens that
$\Gr(A + B_i')= 0$.    This argument shows that if $\Gr_0(A)\ne 0$
and if $A^2 > 0$ then
there is a  decomposition  
$\overline D$ of $K$ with $\Gr_{\overline D} (K)\ne 0$
which contains $A$ as one of its elements.

Now suppose that we are in the situation of (i) and that $A^2 > 0,
\Gr_0(A) \ne 0$.  Then $\Gr_{\overline D}K \ne 0$ and so $\overline D$
must equal $D$, which implies that $A$ must be one of the $B_i$.  To
see that $\Gr(A) = \Gr_0(A)$ we argue by contradiction.  If this is
not true we must have $\Gr_{D''}(A) \ne 0$ for some non trivial
decomposition $D''$ of $A$.  Because (by the compactness theorem)
there are only finitely many classes $B$ with bounded symplectic area
$\om(B)$ which have $J$-holomorphic representatives, we may assume
that $D'' = \{B_1'',\dots, B_p''\}$ consists of indecomposable
elements.  Further, since $A^2 > 0$, one of these elements, say
$B_1''$, must have positive self-intersection number.  Therefore, our
previous argument shows that there is a decomposition ${\overline
D}\,\!''$ of $K$ which contains $B_1''$ such that $\Gr_{{\overline
D}\,\!''}K \ne 0$.  Since $A\ne B_1''$ and $A\cdot B_1'' = (B_1'')^2>
0$ by construction, ${\overline D}\,\!''$ cannot equal ${\overline
D}$: a contradiction.  This proves the first statement in (i).  The
other statements in (i) are now obvious.

Now suppose that there are classes $A_1, A_2$ as in (ii).  We first
claim that $A_i^2 > 0$.  To see this, observe that if $A_1^2 = 0$ for
example, then $K\cdot A_1 = 0$ (since $k(A_1) = 0$.)  But $\Gr(A_2)
\ne 0$ implies $\Gr(K - A_2) \ne 0$, and so, by positivity of
intersections we have $$ K\cdot A_1 = A_2\cdot A_1 + (K-A_2)\cdot A_1
\ge A_2\cdot A_1 > 0, $$ a contradiction.  Hence $A_i^2 > 0$ for both
$i$ and so as above, one can create two different decompositions $D$
of $K$ (one containing each $A_i$) with $\Gr_D(K)$.  The converse is
obvious. \QED

\begin{remark}\label{rem:GW}\rm  If there 
were classes $A_1, A_2$ as in (ii) above, then, by
Ruan--Tian's composition law, there would be a corresponding nonzero
Gromov--Witten invariant in class $A = A_1+A_2$.  Note that $k(A) >
0$.  However, this does not contradict part (i) of Taubes's Structure
Theorem~\ref{Kmin} since the complex structures on the $A$-curves that
we count are not allowed to vary freely but are restricted to be in a
certain cycle in the moduli space arising from the decomposition of
$A$ into $A_1 + A_2$: see~\cite{RT2}.
\end{remark}

\section{Spherical Gromov invariants}

We now develop a theory of \lq\lq spherical" Gromov invariants for
symplectic $4$-manifolds which count the number of ways in which a
class $A$ can be represented by a union of (possibly singular)
$J$-holomorphic spheres.  It is the natural generalization to
disconnected curves of the genus $0$ invariants which were considered
in~\cite{JHOL,RT} and which arise in quantum cohomology.  However,
these spherical invariants are much more limited in scope than the
Gromov invariant considered by Taubes since they vanish on all minimal
symplectic manifolds except for those which are rational and
ruled.\footnote { If the class $A$ represented by the $J$-sphere is
such that $c_1(A) = -K\cdot A > 1$, then this follows by~\cite{IMM},
where it is shown that the only symplectic manifolds that contain such
spheres are blow-ups of rational or ruled manifolds.  Since blowing
down spheres increases $-K\cdot A$, the only case not covered is that
when $M$ is minimal and $-K\cdot A = 1$.  Here we appeal to
Theorem~\ref{Kmin} which shows that $b_2^+ = 1$.  The results of
Liu~\cite{LIU} now show that $M$ is rational or ruled.}  Since the
modified invariant $\Gr'$ (developed in \S3.1 above) is more
appropriate here than the original invariant $\Gr$, we will generalize
$\Gr'$.

Let $(M,\om)$ be a symplectic $4$-manifold.  Given $A\in H_2(M,\Z)$
and $J \in \Jj(M,\om)$, the space $\Mm(A,J)$ of all somewhere
injective $J$-spheres is an oriented manifold of dimension $2(c_1(A) +
2)$.  Hence, if $k = c_1(A) - 1$, the evaluation map $$ e_k: \Mm(A,J)
\times _{{\G}_0} (S^2)^k \to M^k $$ is a map between manifolds of
equal dimension.  It is shown in~\cite{JHOL} that even though the
domain of $e_k$ may not be compact, the map $e_k$ has the structure of
a pseudocycle and hence represents a well-defined element of
$H_{4k}(M^{k})$. (The point is that, by Gromov's compactness theorem,
the image of $e_k$ can be compactified by adding pieces corresponding
to $A$-cusp-curves.  These pieces have to have codimension at least
$2$ and so do not contribute to the homological boundary of $e_k$.)
Moreover, the class $[e_k]$ represented by $e_k$ is independent of the
choice of $J\in \Jj(M,\om)$.  Thus we get a well-defined number by
taking the intersection of $[e_k]$ with a point $(x_1,\dots, x_k) \in
M^k$, or, informally, by counting the number of (unparametrized)
$J$-spheres in class $A$ which go through a fixed generic set of $k$
points in $M$.  This number is the correct number to be called the
{\bf spherical Gromov invariant} \, Gr$_s(A)$ provided that $A$ can be
represented by a somewhere injective $J$-holomorphic immersion of a
$2$-sphere.\footnote{ Recall from~\cite{LOCBE} that any somewhere
injective singular $J$-holomorphic map can be perturbed to a somewhere
injective $J'$-holomorphic immersion for a nearby $J'$.}

Many useful results about $J$-spheres  can be
reformulated in terms of this invariant.  For example:

\begin{prop}\label{prop:sgr1}  If 
 the class $A$ can be represented by a symplectically embedded sphere of
self-intersection number $\ge -1$ then  {\rm Gr}$_s(A)= 1$. 
\end{prop}
\proof{}  The hypothesis implies that $A$ can be represented
by an embedded $J$-sphere for some regular $J$.  By the adjunction
formula, this implies that $c_1(A) = 2 + A\cdot A$.  Then $k = c_1(A)
- 1 > A\cdot A$ and so there can be at most one $J$-sphere through $k$
distinct points.  But there is at least one by hypothesis.  Hence
result.  For more details see~\cite{LM1}.\QED

The following result is proved in~\cite{IMM}.

\begin{prop}\label{prop:sgr2}    If  the moduli space
$\Mm(A,J)$ is non-empty for some regular $J$  then, when $k = c_1(A) - 1\ge 1$,
there are a finite number of $J$ curves through a 
generic set of $k$ points in $M$
and each of these curves contributes $+1$ to  {\rm Gr}$_s(A)$.  In particular,
{\rm Gr}$_s(A)\ge 1$.
\end{prop}

It is not hard to check that the above hypothesis is satisfied
whenever $M$ contains a symplectically immersed $2$-sphere $C$ in
class $A$ with $c_1(C) \ge 2$ whose only singularities are double
points at which the two sheets intersect positively.
\MS

The above definition is fine as far as it goes.  However, as before,
some classes which should have a nontrivial Gromov invariant do not
have connected representatives.  For example, if $A = 2B$ where $B =
[S^2\times \ppt] \in H^2(S^2\times S^2)$, then it follows from
positivity of intersections that, because $B$ always has a
$J$-holomorphic representative and $B\cdot B = 0$, the only $J$-curves
in class $2B$ are $2$-fold coverings of $B$-curves.  Therefore, we
must count curves which may be disconnected.

The only problem in extending our
invariant to this case is to get the correct formula for the number of
points $k$.  Taubes allowed his curves to have arbitrary genus and so used
the number $k(A) = \frac 12 (c_1(A) + A\cdot A)$.  However this is
not appropriate for spheres, since the dimension of the space of
(unparametrized) immersed spheres with $n$ double points is $2n$ less
than the corresponding space of embedded curves of genus $n$.  This is why
we set $k = c_1(A) - 1 $ above.   (The adjunction
formula for connected curves implies that $k(A) - k$ is exactly the genus of
embedded $A$-curves.)  Now the genus of a disjoint union of $p$ spheres is
$-p + 1$  (because the Euler characteristic is $2p$), and so when $p > 1$ we
must change $k$ appropriately.

With this understood, the spherical version of Taubes's definition is as
follows.  

\begin{definition}\label{grs}\rm Given a class $A\in H_2(M)$ consider all pairs
$(k,p)$ where  $$
0 \le k  \le c_1(A)-1, \quad p + k = c_1(A).
$$
  For each such $k$ fix a generic set
$\Om_k$ of $k$ points in $M$ (thus $\Om_0 = \emptyset$) and consider the
set $\Hh_s(A,J)$ of all equivalence classes $(\phi, \Sigma)$  such that

\begin{description}
\item{(a)} $\Si= \cup_i\Si_i$ is a
 Riemann surface which is a disjoint union of 
spheres $\Si_i$;

\item{(b)}
the map $\phi$ 
is $J$-holomorphic on $\Si$ and maps its components to disjoint
curves $\phi(\Si_i)$ in $M$, 
unless $\phi(\Si_i)^2 < 0$ in which
 case the images of components are allowed to
coincide.  Moreover $[\phi_*(\Si)] = A$;

\item{(c)}  elements $\phi,\Si$ and 
$\phi', \Si'$ are equivalent if $\Si = \Si'$ and $\phi'$
is a reparametrization of $\phi$;

\item{(d)}
there is some pair $(k, p)$ as above such that  $\Sigma$ has $p$
components and $\Om_k \subset \phi(\Si)$.    
\end{description}

The  {\bf
spherical Gromov invariant} \, Gr$_s(A)$  is  simply the number of elements in
$\Hh_s(A, J)$.  (Because of
Proposition~\ref{prop:sgr2} we count each pair $(\phi,\Sigma)$ with sign $+1$.)
\end{definition}

\begin{remark}\rm (i) Observe that by condition (b) above exceptional spheres
may occur with multiplicity, so that this is the spherical analog of
the invariant $\Gr'$ rather than $\Gr$.  This multiplicity is handled
somewhat differently than before: we will see below that the
restriction of $\phi$ to each component of $\Si$ has multiplicity one,
however now components may coincide.  To get the spherical analog of
$\Gr$ it is enough to insist that the images of all components are
disjoint.

\NI
(ii) Taubes does not specify in his definition that the images of the
different components of $\Si$ should be disjoint because this is a
consequence of his dimension formula.  However we need to do this
because we have allowed a choice of $(k,p)$.  For example, when $A =
2L$ in $\C P^2$ we wish only to count the unique conic through $5$
generic points, and not the unique pair of lines through $4$ generic
points: see Example~\ref{ex:grs}.
\end{remark}

Part (iii) of the following proposition shows that the new definition
of \,Gr$_s(A)$ agrees with old one (ie the one obtained by counting
{\it connected} curves) in all cases when that was nonzero. Recall
that $\Ee$ denotes the set of classes represented by exceptional
spheres.

\begin{prop}\label{prop:sgr3}  \begin{description}
\item{(i)}  The number  \,{\rm Gr}$_s(A)$ is finite and independent of the
choice of generic $J$ and $\Om$.

\item{(ii)}  If $(\phi, \Si) 
\in \Hh_s(A,J)$, $\phi$ is somewhere injective on every
component of $\Si$.  Moreover, if $\phi(\Si_i)^2 \le 0$,
$\phi|_{\Si_i}$ is an embedding.

\item{(iii)} If the class $A$ can be represented by a somewhere injective
$J$-sphere for some generic $J$, then $\Hh_s(A,J)$ only contains elements with
$p = 1$, and so $\Gr_s(A)$ agrees with the previously defined invariant.
This is the case if $A^2 > 0$ and $E\cdot A \ge 0$ for all $E\in \Ee$. 

\item{(iv)} $\Gr_s(A) = 1$ if all the components of $\Hh_s(A,J)$ are
embedded.
 \end{description}
\end{prop}
\proof{}  Suppose first that $\Si$ has $p$ components and that the
restriction of $\phi$ to $\Si_i$ is an $m_i$-fold covering of a
$B_i$-curve.  Then $A = \sum m_i B_i$, and so there is an associated
element $(\phi', \Si)$ in $\Hh_s(\sum_iB_i, J)$.  Since $\phi'$ is
somewhere injective on $\Si_i$ the curve $\phi_i(\Si_i)$ can go
through at most $k_i = c_1(B_i) - 1$ generic points.  Thus $$ k =
c_1(A) - p = \sum_i k_i = \sum_i c_1(B_i) - p, $$ which is possible
only if all $m_i = 1$.  (Notice that $c_1(B_i) \ge 1$ because, by
assumption, the moduli space $\Mm(B_i,J)/{\G_0}$ of unparametrized
$B_i$-spheres is non empty and so has dimension $2c_1(B_i) - 2 \ge
0$.)  This proves the first statement in (ii).  The second follows in
the usual way fron the adjunction formula and the fact that
$c_1(B_i)\ge 1$.

Suppose now that that the class $A$ may be represented both by a
somewhere injective $A$-sphere and by a union of nonmultiply-covered
spheres in classes $B_1,\dots, B_p$ where $p > 1$, which may coincide
if they are exceptional spheres but otherwise are disjoint.  We first
claim that there cannot be any exceptional spheres among the $B_i$.
For if there were some, in class $E$ say, we would have $E\cdot A < 0$
which contradicts positivity of intersections unless $A = E$.  But
this is impossible because $p > 1$.  Therefore, we may assume that all
the $B_i$ are disjoint.  Therefore, $$ B_i\cdot B_j = 0\; \mbox{ if
}\; i\ne j, \quad A\cdot B_i = B_i\cdot B_i \ge 0\;
\mbox{ for all }\; i.
$$
Suppose further that  $c_1(B_i) > 1$ for some $i$.  Then, by
Proposition~\ref{prop:sgr2}, every generic point of $M$ lies on a $B_i$-curve. 
Therefore, it is possible to have $B_i\cdot B_j = 0$ only if $B_i = B_j$ and
$B_i. B_i = 0$.  By the adjunction formula, $B_i. B_i \ge c_1(B_i) - 2$ with
equality only if $B_i$ is embedded.  It follows that the $B_i$ must be
parallel copies of an embedded curve of self-intersection $0$.  But then $A
=p B_i$ has no $J$-holomorphic somewhere injective  representative.
Therefore, we must have $c_1(B_i) = 1$ for all $i$. 

Since $B_i.B_i - c_1(B_i)$ is even (the adjunction formula again), we
must have $B_i\cdot B_i \ge 1$ so that $A\cdot A \ge p \ge
2$. By~\cite{IMM} this means that $M$ must be a blow up of $\C P^2$ or
$S^2\times S^2$.  The proof of (ii) will be finished by showing that
this manifold does not contain two distinct curves in classes $B_1,
B_2$ satisfying $$ B_1\cdot B_1\ge 1,\quad B_2\cdot B_2\ge 1, \quad
B_1\cdot B_2= 0.  $$ But this follows from the Light Cone
Lemma~\ref{le:light} which holds on all $4$-manifolds with $b_2^+ =
1$.

This proves (iii), and (i) is clear.  Finally (iv) is proved by arguing as in
Proposition~\ref{prop:sgr1}.  Further details are left to the reader.\QED

The following examples illustrate part (iv) of the above proposition.

\begin{example}\rm  (i) Let $X$ be a Riemann surface, and set 
$A = 2F$, where $F$ is the class of the fiber in
$X\times S^2$.  
Then $A$ is represented by a disjoint union of $2$ fibers and
\,Gr$_s(A) = 1$.  

\NI
(ii) Let $M$ be $\C P^2$ blown up at $2$ points, with $L = [\C P^1]$
and $E_1, E_2$ the two exceptional classes.  Then $A = L + E_1 + E_2$
 is represented by the disjoint union of $3$
spheres in classes $L, E_1,$ and $ E_2$ and again it is easy to check that
 \,Gr$_s(A) = 1$  Similarly,  $A = L - E_1 + E_2$
 is represented by the disjoint union of $2$
spheres in classes $L- E_1$ and $E_2$ and \,Gr$_s(A) = 1$.  
\end{example}

Here are some examples which  illustrate the
difference between \,Gr$_s(A)$  and
the Gromov invariant \,Gr$(A)$ (or $\Gr'(A)$).

 \begin{example}\label{ex:grs}\rm  
(i)  If $A = 3L \in H_2(\C P^2)$, then \,Gr$_s(A) =
12$ is  the number of (immersed) $J$-holomorphic rational cubics through
$8$ generic points,  while \,Gr$(A) = 1$ is the number of embedded
$J$-holomorphic  tori through $9$ generic points.   If one blows up a point 
in $\C P^2$ and considers $A = 3L + E$, then we may take $(k,p) = (8,2)$
to obtain \, Gr$_s(A) = $\, Gr$_s(3L)= 12$.  (Note that $A$ itself has no
connected 
$J$-holomorphic representative because $E$ does and $A\cdot E < 0$.)  In
this case, Gr$(A) = \Gr'(A) = 1$ is represented by 
the union of the unique torus
through $9$ generic points with the exceptional curve.  \MS

\NI
(ii)  If $A$
is the class of $T^2\times \ppt$ in $T^2\times S^2$ then \, Gr$_s(A)$ is
obviously $0$ since $A$ has no spherical representatives, while \, Gr$(A) =
\Gr'(A) = 
2$: see Lecture~5.
\end{example}

However it is easy to check that the two invariants do agree in the following
situation.

\begin{prop}  If the class $A$ can be
represented by a disjoint union of embedded $J$-spheres then $\Gr_s(A)=
\Gr'(A)$.  This is the case whenever $\Gr'(A)$  is
calculated using spheres. 
 In particular  $\Gr_s(A) = \Gr'(A) (= \Gr(A))$ when $A$ is the 
class of the fiber in a ruled symplectic manifold or when $A = [\C P^1] \in
H_2(\C P^2)$. 
\end{prop}
\proof{}  Suppose that $A$ can be
represented by a disjoint union of embedded $J$-spheres.
  If all these components
are exceptional spheres then $A$ is a sum of elements from $\Ee$ and the result
follows from Proposition~\ref{goodc2}.  Otherwise, 
$M$ must be a blow-up of a rational or ruled mnifold and the result 
follows from Propositions~\ref{prop:b+1} and~\ref{prop:sgr3}.\QED

\begin{remark}\rm  According to the perspective of Ruan--Tian in~\cite{RT2}, 
instead of counting immersed $J$-spheres with $d$ double points going
through a set of $k$ points, one can resolve the singularities and
count instead the number of genus $d$ curves through $k$ points whose
complex structures are constrained to lie in an appropriate cycle in
the moduli space (of complex structures on genus $d$ curves).  Thus
the spherical invariant $\Gr_s$ is really a special case of their
Gromov--Witten invariants.  There are many open questions about these
invariants.  For example, given a minimal symplectic manifold with
$b_2^+ > 1$ what is the minimal genus of a curve for which some
Gromov--Witten invariant is nonzero?  It is not even known whether
there are non trivial examples of these invariants on manifolds with
$b_2^+ > 1$, i.e.  ones which cannot be reduced to Taubes's Gromov
invariants: cf Remark~\ref{rem:GW}.
\end{remark}

\section{Calculating Gromov invariants of  tori}

We begin with an example illustrating what can happen with tori, and
then outline Taubes's method for counting them.  Finally we discuss
the Gromov invariant of the fiber class in an elliptic surface without
multiple fibers, and show how to calculate it using a sum formula.  A
good general reference for facts about complex surfaces is Griffiths
and Harris~\cite{GH}.

\subsection{Tori in $S^2\times T^2$}\label{ss:tor}

Let $M = S^{2}\times T^{2}, $ and $B={\rm pt}\times T^{2}.$
Suppose that $B$  is  represented by an embedded torus. Since 
$$
\dim
\Mm (B, J, 1) = 2 (c_{1}(B) + g - 1) + \dim G_{1} = \dim G_{1},
$$ 
the dimension of the unparametrized moduli space is $0$.  In other
words, regular tori in class $B$ are isolated.  Thus the product
complex structure on $T^2\times S^2$ is not regular. To find a regular
$J$.  realise $S^{2}\times T^{2}$ as the projectivization $ {\bf P} (L
\oplus \C)$ of the rank-$2$ bundle $L\oplus\C$, where $L \to T^{2}$ is
a nontrivial holomorphic line bundle with $c_{1} = 0.$ We claim that
with this complex structure $J_L$ the manifold $M$ has exactly $2$
$J_L$-tori in class $B$ which both count with $+1$.  Hence $\Gr(M,B) =
2$.

To see this, observe first that if $L$ had a holomorphic section, this
section cannot vanish anywhere, since every intersection with the zero
section counts positively by positivity of intersections.  Hence our
$L$, which is nontrivial by assumption, does not admit nonzero
holomorphic sections.  Moreover, the only holomorphic sections of
${\bf P} (L \oplus \C)$ are $\left[L\oplus 0\right],$ the section at
\lq\lq infinity'' and $\left[ 0 \oplus \C\right]$, the \lq\lq zero
section''.  (To see this, observe that such sections are in bijective
correspondence with line subbundles $E$ of $L\oplus \C$.  But if $E
\ne L\oplus
\{0\} $ or $\{0\}\oplus \C$, $E$ gives rise to a nontrivial homomorphism from
$L^*$ to $\C$, which does not exist.  For more details on this kind of
argument see~\cite{RULED}.)  Further it is easy to check that the
normal bundle of $\left[L\oplus 0\right]$ is isomorphic to $L^*$,
while that of $\left[ 0 \oplus \C\right]$ is isomorphic to $L$.  This
is obvious for the section $[0\oplus\C]$, and follows for
$\left[L\oplus 0\right]$ since the latter can also be identified with
$\left[\C\oplus 0\right]$ in $$ {\bf P}\, (\C
\oplus L^*) = 
 {\bf P} \,\left(L^*\otimes (L \oplus \C)\right) = {\bf P} (L \oplus
\C).  $$ Since these normal bundles $\nu$ are nontrivial,
$H^{1}_{\overline{\partial}}(T^{2}, \nu) = 0$ for both sections, which
implies that these curves are regular.  Hence $J_L$ is regular for the
class $B$.  Moreover, since $J_L$ is integrable both tori count with a
$+$ sign.  Hence $\Gr(B) = 2$.

For generic $L$, I next claim that there are no embedded
$J_L$-holomorphic tori in the class $2B$.  For if there were, there
would be a double cover map $\psi:T^2\to T^2$ such that this torus
pulls back to a torus $T'$ say in class $\ppt\times [T^2]$ in the
manifold $(S^2\times T^2, J_{L'})$ where $L' = \phi^*(L)$.  Since
$L'\ne \C$ for generic $L$ and $T'$ is neither the section at zero nor
that at infinity, this is impossible.  However, there are three
representatives $(\phi,\Si)$ of the class $2B$, namely double covers
of each of the $B$-curves and a disconnected curve with $2$ distinct
components.  It is not hard to check that these are regular (for
generic $L$).  Again, adopting the principle that (regular)
holomorphic objects always count with $+1$\footnote { For
nonmultiply-covered curves, this principle is justified by the fact
that in the case all evaluation maps are holomorphic maps between
complex manifolds and so all intersection numbers are positive.}, we
find that $\Gr(2B) = 3$.  More generally, we have

\begin{lemma}  If $B = [\ppt\times T^2]\in H_2(S^2\times T^2)$ then $\Gr(kB)
=  k+1$ for $k > 0$.
\end{lemma}
\proof{}(Sketch)   Arguing as above 
we see that for generic $L$ the only connected
curves in class $pB$ are $p$-fold covers of the sections at zero and
infinity.  Hence there are exactly $k+1$ ways of representing the
class $kB$, and each counts with a $+1$. \QED

\NI
{\bf Note:}
The above result can be fully justified using Taubes's work in~\cite{TAUTOR}.
It  is also compatible with the calculation via Seiberg--Witten invariants:
see Li-Liu~\cite{LL1,LL2}.\MS

Now Taubes shows in~\cite{TAUTOR} how to define an almost complex
structure $J_1$ on $S^2\times T^2$ which is also regular for the class
$B$ but which admits $4$ $J$-holomorphic tori in class $B$: the two
above plus a cancelling pair, one which occurs with a $+$ sign and one
with a $-$ sign.  Moreover he shows that this $J_1$ admits no embedded
tori in class $2B$.  (Such an example was worked out independently by
Lorek in~\cite{LOR}.)  Suppose that the correct way to count
multiply-covered tori is simply to assign a $\pm 1$ to each multiple
covering according to some rule and otherwise to follow the scheme
laid out in Definition~\ref{def:gr}.  Then of the $6$ disconnected
representatives of $2B$ three occur with a $+1$ and three with a $-1$
and so the net contribution is $0$.  But there are $4$ doubly covered
curves, and there is {\it no} way to assign the numbers $\pm 1$ to
these four curves to make them give $3$.

\subsection{Taubes's method for counting tori}

Looking at examples like this, Taubes realised that to take proper
account of the way in which a multiply covered torus contributes to
the Gromov invariant one has to look at more than just the orientation
of the underlying embedded curve.  To understand why this is so,
consider a generic path $J_t, 0\le t\le 1,$ of $\om$-tame almost
complex structures.  Any regular embedded $J_0$-holomorphic torus
$C_0$ in class $A$ is the endpoint of a path $C_t$ of
$J_t$-holomorphic tori in class $A$.  As $t$ increases, two kinds of
bifurcations occur.  One is the birth or death of a pair of tori (one
$+$ and one $-$ as above).  The other is more complicated: a torus
$T_t$ in class $2A$ can split off from the basic path $C_t$.  (A
beautiful explicit model for this bifurcation is described by Lorek in
\S4 of~\cite{LOR2} which exhibits it as a quotient of a standard birth
bifurcation.)  The double covering map $T_t \to C_t$ is classified by
one of the $4$ elements of $\tau\in H_1(C_t,\Z/2\Z)$.  Taubes's basic
insight is to realise that in order to determine the weight to assign
to a $k$-fold cover of the torus $C_0$ one must keep track of the sign
of the determinants ${\rm Det}\,D_\tau$ of the linearized
Cauchy-Riemann operators $D_\tau$ on the normal bundle $\nu_{C_0}$
twisted by all {\it four} of the elements $\tau$.  (Taubes defines
these signs in terms of a suitable spectral flow.)  When $\tau = 0$
this sign corresponds exactly to the sign $\pm 1$ that we were using
above to weight a curve, but when $\tau\ne 0$ this is new data.  Thus
each regular curve gets one of eight possible labels $$ (\pm, i),
\quad i= 0,1,2,3 $$ where the label $(\pm,i)$ means that ${\rm
Det}\,D_\tau = \pm 1$ when $\tau = 0$ and that exactly $i$ of the
other three signs are equal to $-1$.  Observe that if $J$ is
integrable near the curve $C$ its labels are always $(+,0)$; in other
words, all determinants are $+1$ in accordance with our previous
positivity principle.

Taubes establishes:
\SS

\NI
{\bf Birth rule}:\, when two curves 
are born the pair has labels $(\pm, i)$ for some
$i$.    
\SS

\NI{\bf Bifurcation rules}:\,
these describe how the labels for $C_t$ change as $C_t$ passes through
a bifurcation in which a torus $T_t$ corresponding to $\tau$ splits
off, and they also give a formula for the label of the new torus
$T_t$.
\SS

Using this information, he proves that the following method for
weighting tori gives rise to a way of counting tori which is invariant
under symplectic deformation.  For each label $(\pm,i)$, let $$
f_{(\pm,i)}(t) = \sum_{k\ge 0} n_kt^k $$ be the generating function
for its contribution to the Gromov invariant.  This means that if a
torus in class $A$ has label $(\pm ,i)$ then its contribution to the
count of tori in class $kA$ is $n_k$.  Then
\begin{eqnarray*}
f_{(-,i)} & = & \frac 1 {f_{(+,i)}},\\\vspace{.1in}
f_{(+,0)} = \frac 1 {1-t},& &
f_{(+,1)}= 1+t,\\\vspace{.1in}
f_{(+,2)}= \frac{1+t}{1+t^2}, & &
f_{(+,3)}= \frac{(1+t)(1-t^2)}{1+t^2}.
\end{eqnarray*}

We now  complete the definition of $\Gr_0(A)$ that was begun in
Definition~\ref{def:gr1}. 

\begin{definition}[Taubes]\label{def:gr2}\rm  Let  $A\in H_2(M,\Z)$  be a
primitive element such that $A^2 = 0$ and $K\cdot A = 0$, and , for
generic $J$, let $\Tt(A,m)$ denote the set of embedded, regular
$J$-holomorphic tori in class $mA$.  Then the Gromov invariant
$\Gr(M,kA) = \Gr(kA)$ for $kA$ in $M$ is the coefficient of $t^k$ in
the power series $$
\prod_{m\le k} \left(\sum_{C\in \Tt(A,m)} f_{\ell_C}(t^{m})\right),
$$ where $\ell_C$ is the label for $C$.  Moreover, we set $\Gr_0(kA) =
\Gr(kA)$ in this case. \end{definition}

To illustrate this, consider the extra $\pm$ pair of $J_1$-holomorphic
tori which cause trouble in the example above.  In the construction,
it turns out that the $+$ torus has type $(+,0)$, so that the $-$
torus has type $(-,0)$.  Hence all $k$-fold coverings of the $+$ torus
count with $+1$, while $k$-fold coverings of the $-$ torus with $k >
1$ contribute $0$ (because the above rules give $f_{(-,0)} = 1-t$.)
To see how this gives rise to an invariant which is independent of
$J$, consider a generic path $J_t$ from an integrable structure $J_0$
on $T^2\times S^2$ to $J_1$, and suppose, for simplicity, that there
is a single bifurcation point on this path at which the $\pm$ pair
$C_{\pm}$ is created.  Then, the bifurcation creates $6$ new elements
that contribute to $\Gr(2A)$.  These are the disjoint union of one of
the new tori with one of the two old ones, together with two elements
involving just the new tori $C_\pm$, namely the double cover of $C_+$
and the disjoint union of the pair of them. It is now easy to check
that the net contribution to $\Gr(2A)$ of these $6$ new elements is
zero.  Observe, in particular, that it is impossible to get a well
defined, consistent invariant, which counts only connected curves in
class $2A$.  This is why we lumped all parallel tori together when
defining the decomposition of $A$ in Definition~\ref{def:decomp}.

\subsection{Elliptic surfaces}

Consider the rational elliptic surface $V(1) = \C P^{2} \# 9
\overline{\C P}^2.$ This may be understood as follows.  Take nine
points in general position in $\C P^2$. There exists a unique torus (a
cubic curve) through those nine points in the class $3 L = 3 \left[\C
P^{1} \right].$ (One can check this simply by looking at the
corresponding system of $9$ equations in the $9$ unknown coefficients
of the cubic.)  Blow up the nine points to obtain a torus $T$ in class
$F = 3 L - E_{1} -\cdots - E_{9},$ with zero self-intersection $F\cdot
F = 0.$ Since all elliptic surfaces do embed in $\C P^2$ we may
suppose that the induced complex structure $j$ on this torus $T$ is
generic.  Then, if we choose the $9$ points on this torus $T$ so that
they are also generic, it is not hard to check that $T$ is regular, ie
that its holomorphic normal bundle satisfies $H^1(T, \nu) = 0$.  This
shows that the Gromov invariant $\GR (F) = 1.$ In fact, it follows
from Seiberg--Witten theory that $\Gr(B) = 1$ for every class $B = nL
-\sum_im_iE_i$ such that $n>0$ and $B^2\ge 0$, as well as for the
classes $B = E_i$.

Now choose nine points which lie on two cubic curves $ C_{i} = \{
f_{i} = 0\},$ for $ i =1,2.  $ Then in fact the points lie on the one
parameter family of cubics $$ C_\la = \{f_1 + \la f_2 = 0\}.  $$
Moreover, each point of $\C P^2$ except for the $9$ points of
intersection lies on exactly one of these cubics (provided that we
allow $\la = \infty$ to include the cubic $C_2$).  Thus there is a
well-defined map of the complement of the $9$ points to the parameter
space $\C P^1$.  It is not hard to check that it extends to a smooth
holomorphic map $V(1)\to \C P^1$.  If $C_1$ and $C_2$ are generic,
this is a singular fibration in the sense that the fiber over all but
finitely many of the points $\la$ is the nonsingular cubic $C_\la$.

A manifold which fibers like this over $\C P^1$ is called an elliptic
surface.  $V(1)$ is the simplest elliptic surface.  We can construct
others by Gompf's construction of the fiber connected sum.  For
example $V(2) $ is obtained from two copies of $V(1)$ by removing a
neighborhood of a fiber in each copy and then gluing the boundaries by
a suitable orientation reversing symplectomorphism.  To be precise,
recall that by the symplectic neighborhood theorem a small open
neighbourhood $N$ of a fiber $T^{2}$ in $V(1)$ is symplectomorphic to
a product $T^{2}\times D^{2}$.  By making $N$ smaller if necessary, we
may assume that this product structure extends to a neighborhood $W$
of the boundary of $V(1) - N$.  Thus $W$ is symplectomorphic to
$T^2\times A$ where $A\subset\R^2$ is an annulus.  The surface $V(2)$
is then defined by $$ V(2) = (V(1) - N)\cup _{W = W'} (V(1) - N'), $$
where the identification of $W $ with $W'$ is via a symplectomorphism
of the form $$ {\rm id}\times \phi: T^2\times A \to T^2\times A'.  $$
Here $\phi: A\to A'$ is area preserving and turns the annulus $A$
inside out so that it maps the boundary of $V(1) - N$ into the
interior of $V(1) - N'$.  When repeated this construction yields a
family of elliptic surfaces $V(n) = V(n-1) \# _{T} V(1).$

The construction given above takes place in the symplectic rather than
holomorphic category.  However it is not hard to construct $V(n)$ (eg
as a branched cover) in a way that makes clear that it does have a
complex structure.  Interestingly enough, when $n > 1$ {\it all}
complex structures on $V(n)$ give it the structure of an elliptic
surface.  (This is {\it not} true when $n = 1$.)  When $n = 2$ we get
a $K3$ surface.  This surface is analogous to the $4$-torus $T^4$ in
that generic complex structures on it admit no holomorphic curves at
all.  Therefore all its Gromov invariants vanish.  Moreover, as in the
case $n = 1$, a generic complex structure on $V(2)$ is regular in the
Fredholm sense, ie the moduli spaces of holomorphic curves are
manifolds of the right dimension and so can be used to calculate
Gromov invariants.

However, when $n > 2$ no complex structure on $V(n)$ is regular in the
sense of Fredholm theory.  For, regular $J$-holomorphic tori of zero
self-intersection are isolated.  However, because the moduli space of
holomorphic fibers in $V(n)$ is a manifold (albeit of too high a
dimension) one can still try to use it to calculate the Gromov
invariant.  In fact, Ruan shows in~\cite{RU} that the contribution of
a {\it compact} component of the moduli space to the Gromov invariant
of a class with formal dimension (or index) $0$ is precisely its
oriented Euler characteristic.  One cannot simply apply that result
here though, since the moduli space of holomorphic tori in the class
of the fiber is not compact.  (Its ends are the singular fibers.)
Nevertheless, as Ruan pointed out to me, there is a natural
holomorphic compactification of this moduli space which has Euler
characteristic $2-n$.  Thus, this heuristic calculation suggests that
$\Gr(V(n), F)$ should be $2-n$.  Note that this is a negative number.
However, this does not contradict the positivity principle that we
were using before because the holomorphic objects here are not regular
in the Fredholm sense.  In fact, one can think that this is a reason
why no complex structure on $V(n)$ can be regular when $n > 2$: the
only way that a negative number can be holomorphically represented is
through a nonregular family with too large dimension and negative
global twisting.

\subsection{The Gromov invariants of a fiber sum}
 
One can show that $\Gr(V(n), F) = n-2$ from Seiberg--Witten theory,
using the fact that $V(n)$ is K\"ahler.  Here we show how to calculate
it directly from the construction of $V(n)$ as a fiber sum.  The full
details, together with a general statement valid for any connected sum
along tori, have been worked out by Lorek~\cite{LOR2}.  I wish to
thank him for useful discussions and in particular for pointing out
the role of the so-called boundary classes.

The basic idea is to consider Gromov invariants for compact manifolds
$(X,\om)$ which have boundary components diffeomorphic to $S^1\times
T^2$.  We will always suppose that the symplectic form restricts on
the boundary to the pullback of the usual area form on the torus
$T^2$.  Hence, by the symplectic neighborhood theorem, we may identify
a neighborhood of each boundary component with $$ (N\times T^2, \om) =
([-1,0]\times S^1\times T^2, du\wedge d\theta + ds\wedge dt), $$ where
the boundary is at $u = 0$.  We will take $A = \ppt\times T^2$.  In
order for the Gromov invariant to be well-defined we must make sure
that $J$-holomorphic tori cannot escape through the boundary of $X$.
Therefore we will only consider almost complex structures $J$ on $X$
which have the following standard form on a neighborhood of each
boundary component of $X$.

Given a function  $\be:[-1,0]\to \R$ which satisfies the following
conditions: \SS

\NI
(i)    $0 <|\be(u)| < 1$ for  $u\in [-1,0)$ and $\be(0) = 0$;

\NI
(ii)  $\be$ has isolated critical points, and all
 its derivatives vanish at $u = 0$;
\SS

\NI
we define $J_\beta$ by setting:
\begin{eqnarray*}
\Jb (\p_u) = \p_\theta,& &\Jb (\p_\th) = -\p_u,\\
\Jb(\p_s) = \p_t +\be(u)\p_\th,& & \Jb(\p_t) = -\p_s +\be(u)\p_u.
\end{eqnarray*}
It is easy to check that condition (i) implies that $\Jb$ is
$\om$-tame.  Observe also that for those $c$ with $\be(c)$ rational,
the $3$-torus $\{u = c\}$ is foliated by $\Jb$-holomorphic tori of
\lq\lq slope" $\be(c)$.  More precisely, this foliation is the kernel
of the $1$-forms $$ du = 0, \quad d\th + \be(c)dt = 0, $$

\begin{definition}\rm  An almost complex structure $J$  on $(X, \om)$ is said to be
{\bf compatible with the boundary} of $X$ if each boundary component
of $X$ has a neighbourhood $N$ such that the triple $(N,\om, J)$ is
symplectomorphic to $([-1,0]\times S^1\times T^2, du\wedge d\theta +
ds\wedge dt, \Jb)$ for some $\Jb$ which satisfies the above
conditions.  Elements $B\in H_2(X,\Z)$ which are in the image of some
inclusion map $H_2(N)\to H_2(X)$ will be called {\bf boundary classes}
of $X$. \end{definition}

The next lemma shows that the Gromov invariant $\Gr(X,A)$ is
well-defined when $A\in H_2(X)$ is any class that has zero
intersection with all boundary classes $B$.

\begin{lemma}\label{bound}  Suppose that $X $  has
boundary components $S^1\times T^2$ as above, and suppose that $J_X$
is a generic $\om$-tame almost complex structure on $X$ which is
compatible with the boundary of $X$.  Then, for any class $A$ such
that $A\cdot B = 0$ for all boundary classes $B$, the number of
$J$-holomorphic tori in class $A$ is independent of the choice of $J$.
\end{lemma}
\proof{}  We suppose for simplicity that
 $X$ has a single boundary component.  The
proof in the general case is similar.

Let $F = \ppt\times T^2$.  We first show that $\Gr(X,[F])$ is
well-defined.  Observe first that the boundary $3$-torus $\{u = 0\}$
is foliated by $J_X$-holomorphic tori in class $[F]$, and so there is
a corresponding circle in the moduli space $\Mm(J_X,[F])/G$.  By
Ruan's result in~\cite{RU} the contribution of this circle to the
Gromov invariant is its Euler characteristic, namely $0$.  Because
$|\be(u)| < 1$ none of the other compact leaves in the $3$-tori $u =
c$ can represent $[F]$.  The next important fact is that no other
$J$-holomorphic torus can intersect the boundary region $N\times T^2$.
For if it did, it would have to cross one of the $3$-tori $\{u = c\}$
with $\be(c)$ rational, and hence it would have to intersect one of
the compact leaves of the $J$-holomorphic foliation of this $3$-torus.
Since these leaves lie in a boundary class $B$, this contradicts the
fact that $F\cdot B = 0$.  Therefore any other $J_X$-holomorphic torus
in class $[F]$ has to be contained in the complement of $N\times T^2$.
This means that the boundary region $N\times T^2$ functions somewhat
like a pseudo-convex boundary, containing the $J$-holomorphic curves
in class $[F]$.  In particular, the moduli space of $[F]$-curves in
$X$ is compact.  Moreover, it follows from the usual transversality
arguments that we can therefore calculate $\Gr(X,[F])$ by counting the
elements in the moduli space of $J$-holomorphic $[F]$-tori for a
generic element $J$ of the set $$
\Jj_N = \{ \om\mbox{-tame}\; 
J \;{\rm on}\; X: J = J_X \;{\rm on }\; N\times T^2\}.
$$ 
(For example, if you look at the proof of Proposition 3.4.1
in~\cite{JHOL}, you see that in order to prove that the universal
moduli space $\Mm(A,\Jj_N)$ is a manifold it suffices to be able to
make $J\in \Jj_N$ generic {\em somewhere} on the image of every curve
in $\Mm(A,\Jj_N)$, but not necessarily everywhere on $X$.)

This proves the result when $A$ is the fiber class $ [F]$.  The
argument for other $A$ is similar.\QED

\begin{prop}\label{ellip} 

\begin{description}\item{(i)}  
$\Gr(D^2\times T^2, [F]) = 1$ where $F = \ppt\times
T^2$.

\item{(ii)}   If $N(F)$ is a 
neighborhood of the fiber $F$ in the elliptic surface
$V(n)$ then
$$
\Gr(V(n) - {\rm Int\,} N(F), [F]) = 1-n,\quad \Gr(V(n),[F]) = 2-n.
$$ 
\end{description}
\end{prop}
\proof{}  In~\cite{LOR2} Lorek calculates  $\Gr(D^2\times T^2, [F])$ 
by explicit construction of a suitable $J$ for which one can count the
tori.  We will give a nonexplicit proof which uses the fact that we
know that $\Gr(S^2\times T^2, [F]) = 2$.

 First, let consider the situation when $F$ is a symplectic torus with
$F^2 = 0$ in a closed manifold $Y$.  Then $F$ has a neighborhood
$N(F)$ which is symplectomorphic to the product $D^2\times F$.
Consider the decomposition of $Y$ into $$ Y = Y_0 \cup N(F),
\quad\mbox{where}\;\; Y_0 = Y - {\rm Int}\,N(F), $$ and let $J$ be a
generic almost complex structure on $Y$ formed by putting together
almost complex structures on $Y_0$ and on $N(F)$ which are compatible
with their boundaries.  (Observe that $J$ is smooth because of
condition (ii) on $\be$.)  Then, the moduli space of $J$-holomorphic
tori in class $A = [F]$ splits into three parts: the circle of tori
along the boundary, the set $\Mm_Y$ of tori in $Y_0$ and the set
$\Mm_N$ of tori in $N(F)$. Clearly, $$
\Gr(Y,[F]) = \#\Mm_Y \,+ \,\# \Mm_N,
$$
where one counts the number $\#$ with appropriate signs.  

One complicating factor that we have to consider here is that the
inclusion $Y_0 \to Y$ need not induce an injection on $H_2$. For
example, if we take $Y = V(1)$, the torus in the boundary $S^1\times
F$ which is the kernel of the $1$-form $d\th + dt = 0$ is homologous
to $F$ in $Y$ but not in $Y_0$.  Hence the set of tori in $\Mm_Y$ need
not all lie in class $A$.  Thus, although $\# \Mm_N = \Gr(N(F),[F])$,
the number $\#\Mm_Y$ is, in general, a sum of Gromov invariants.
However, if $Y = S^2\times T^2$, then $Y_0 = D^2\times T^2 = N(F)$ and
this problem does not arise.  Thus we find $$ 2 = \Gr(S^2\times T^2,
[F]) = 2 \Gr(D^2\times T^2,[F]).  $$ This proves (i).

We prove (ii) by induction.  Above we showed that there was a regular
complex structure on $V(1)$ which had exactly one holomorphic torus in
class $[F]$.  Hence $\Gr(V(1),[F]) = 1$.  Decompose $Y = V(1)$ into
$Y_0\cup N(F)$ as above, and consider the set of tori $\Mm_Y$.  It
would theoretically be possible that these tori would give rise to
nonzero Gromov invariants $\Gr(Y_0, A)$ where $A$ is some boundary
class, since $\Mm_Y$ might contain a plus torus in one class $B_+$ and
a minus torus in another class $B_-$ which would cancel out in $Y$ but
not in $Y_0$.  However, if we double $Y_0$ we get the $K3$ surface
$V(2)$ and we know that all invariants vanish for that.  Since the
inclusion of $Y_0$ into its double induces an injection on homology,
we know that the invariants for the double are exactly twice those for
$Y_0$.  Hence the invariants for $Y_0$ must vanish.  This shows that
$$
\Gr(V(1) - {\rm Int\,} N(F), B) = 0,
$$ 
for all boundary classes $B$ including $B = [F]$.

Thus (ii) holds for $n = 1$.  Here is the rough idea for the inductive
step.  Think of $V(n+1)$ as the fiber sum $V(n)\# V(1)$.  Thus we get
$V(n+1)$ by cutting out a neighborhood of a fiber in each of $V(n)$
and $V(1)$ and gluing the remaining pieces together.  We saw above
that if we we make $J$ compatible with the boundaries of the pieces
then the neighborhoods of the cut-out fibers each contain exactly one
$+$ torus.  As we saw above, this uses up the $+$ torus in $V(1)$
leaving it with no tori.  Now $V(n)$ may not have any $+$ tori, and so
to do the cutting we must create a $\pm$ pair (without introducing any
new tori) and then cut out the $+$ one.  This creates an extra $-$
tori in $V(n+1)$, which gives the result.

To make this precise we just have to see that it is possible to make a
 new $\pm$ pair of tori in the fiber class $[F]$ without creating any
 new tori in boundary classes.  Inductively, we can suppose that $J_n$
 is an almost complex structure on $V(n)$ which has the form $J_\be$
 on some set symplectomorphic to $N = [0,1]\times S^1\times T^2$, and
 which is generic outside $N$.  (This means that $J_n$ can be used to
 calculate the Gromov invariant of $V(n)$.)  It follows from
 Lemma~\ref{bound} that $\Gr(N,B)= 0$ for all boundary classes $B$.
 Let $J_n'$ be an almost complex structure which equals $J_n$ outside
 a compact set in ${\rm Int\,} N$ and is such that it is compatible
 with the boundary of some subset $P$ of \,${\rm Int\,} N$
 symplectomorphic to $D^2\times T^2$.  Then, we claim that $$
\Gr(N - {\rm Int\,} P, [F]) = -1,\quad \Gr(N - {\rm Int\,} P, B) = 0,
$$ 
where $B$ is any boundary class other than $[F]$.  The first
statement holds because $\Gr(P,[F]) = 1$, and the second holds
because, as before, there would otherwise be nontrivial invariants for
the $K3$ surface.  Hence, when forming $V(n+1)$ we can cut out $P$ and
glue in $V(1) - N(F)$, a process which leaves us with one more $-$
torus than there was in $V(n)$.\QED



\begin{thebibliography}{9999}

\bibitem{FS}  R. Finushel and R. Stern, 
Rational blow-downs of smooth $4$-manifolds,
UCIrvine preprint (1995)

\bibitem{GOM}  R. Gompf:  A new construction for symplectic
      $4$-manifolds,  {\it Annals of Mathematics}, {\bf 142} (1995),
       527--595.

\bibitem{GH}  P.~Griffiths and J.~Harris,   
    {\it Principles of Algebraic Geometry},
    Wiley, New York, 1978. 


\bibitem{LM}  F. Lalonde and D McDuff, The classification of ruled
       symplectic $4$-manifolds, Stony Brook preprint (1995).

\bibitem{LM1}  F. Lalonde and D McDuff, $J$\/-curves and the 
       classification of rational and ruled symplectic
        $4$\/-manifolds, preprint (1995), to appear in the Proceedings
        of the 1994 Newton Institute Symplectic Geometry conference,
        ed Donaldson and Thomas, Camb Univ Press (1996).

\bibitem{LL1} T.J. Li and A. Liu, General wall 
	crossing formula, Harvard preprint (1995).

\bibitem{LL2} T.J. Li and A. Liu, Symplectic structure on ruled surfaces and
    generalized adjunction formula, Harvard  preprint (1995).

\bibitem{LIU} A. Liu, Some new applications of general wall crossing
   formula, Harvard  preprint (1995).

\bibitem{LOR}  W. Lorek, Nongeneric almost complex structures in symplectic
$4$-manifolds, Stony Brook preprint (1995)

\bibitem{LOR2}  W. Lorek,  A note on Gromov invariants of connected sums of
symplectic manifolds, Stony Brook preprint (1996)

\bibitem{RR} D.~McDuff, 
      The structure of rational and ruled symplectic $4$-manifolds,
      {\it Journ. Amer. Math. Soc.} {\bf 3} (1990), 679--712;
      Erratum: {\it Journ. Amer. Math. Soc} {\bf 5} (1992), 987--988.

\bibitem{LOCBE} 
D. McDuff,   The local behaviour of holomorphic curves in almost complex
$4$-manifolds. {\it Journal of Differential Geometry}, {\bf 34}, (1991) 311--358. 

\bibitem{IMM} D.~McDuff, 
    Immersed spheres in symplectic $4$-manifolds,
    {\it Ann. Inst. Fourier, Grenoble} {\bf 42} (1992), 369--392.

\bibitem{RULED} 
D. McDuff,  Notes on ruled symplectic $4$-manifolds. 
{\it Transactions of the American Mathematical Society\/} {\bf 345} (1994)
623--639.

 
\bibitem{JHOL} 
D. McDuff and D.A.Salamon,   {\it $J$-holomorphic curves 
and quantum cohomology}. University Lecture Series, American Mathematical Society,
Providence, RI. (1994).

\bibitem{INTRO} 
D. McDuff  and D.A.Salamon, {\it Introduction to Symplectic Topology}, OUP (1995).

\bibitem{MS} 
D. McDuff and D.A.Salamon,  A survey of symplectic $4$-manifolds with $b_2^+ =
1$,  Proceedings of the 1995 G\"okova conference, ed Akbulut and Onder, International
Press (1996).

\bibitem{OO} H. Ohta and K. Ono, Note on symplectic $4$-manifolds with
      $b_2^+ = 1$, II, preprint (1995).

\bibitem{RU}  Y. Ruan, Symplectic Topology and Extremal Rays, {\it Geom. 
Funct. Anal.} {\bf 3}, (1993) 395--430.

\bibitem{RUAN2} Y. Ruan, Symplectic topology and Complex Surfaces,
{\it Geometry and Topology on Complex Surfaces}, ed Mabuchi, Noguchi, Ochial,
World Scientific Publications, Singapore (1994).

\bibitem{RT} Y. Ruan and G. Tian, A mathematical theory of Quantum
Cohomology, {\it J. Diff. Geom} (1995)

\bibitem{RT2} Y. Ruan and G. Tian, Higher Genus Symplectic Invariants and
Sigma model coupled with gravity, MIT preprint 1996.

\bibitem{SAL}  D. Salamon, {\it Spin Geometry and Seiberg--Witten
invariants}, in preparation.

\bibitem{TAU0}  C. H.~Taubes,  The Seiberg--Witten invariants and
      symplectic forms, {\it Math. Res. Letters }{\bf 1} (1994), 809--822

\bibitem{TAU00}  C. H.~Taubes,  More constraints on symplectic forms from
       Seiberg--Witten invariants, {\it Math. Res. Letters
       }{\bf 2} (1994), 9--14.

\bibitem{TAU}  C. H.~Taubes, The Seiberg--Witten and the Gromov invariants,
      Harvard preprint (1995).

\bibitem{TAU1}  C. H.~Taubes, From the Seiberg--Witten equations to
      pseudo-holomorphic curves, Harvard preprint (1995).

\bibitem{TAUTOR} C. H.~Taubes, Counting pseudoholomorphic submanifolds in
dimension $4$, Harvard preprint (1995).

\bibitem{WIT} E. Witten,   Monopoles and $4$-manifolds, preprint,
hep-th/9411102,  (1994).


 \end{thebibliography}
\end{document}